\documentclass[aihp]{imsart}

\RequirePackage{amsthm,amsmath,amsfonts,amssymb}
\RequirePackage[numbers,sort&compress]{natbib}
\RequirePackage[colorlinks,citecolor=blue,urlcolor=blue]{hyperref}
\RequirePackage{graphicx}
\usepackage{indentfirst}
\usepackage{caption}
\usepackage{mathrsfs}

\startlocaldefs
\theoremstyle{plain}

\newtheorem{theorem}{Theorem}[section]
\newtheorem*{theorem 2.12}{Theorem 2.12}
\newtheorem{lemma}[theorem]{Lemma}
\newtheorem{definition}[theorem]{Definition}

\newtheorem{corollary}[theorem]{Corollary}
\theoremstyle{remark}
\newtheorem{remark}{Remark}

\endlocaldefs

\begin{document}
\begin{frontmatter}
\title{Persistence of Invariant Tori for Stochastic Nonlinear Schr{\"{o}}dinger in the Sense of Most Probable Paths}
\runtitle{Persistence of Invariant Tori for Stochastic Nonlinear Schr{\"{o}}dinger}

\begin{aug}
	\author[A]{\inits{F.}\fnms{Xinze}~\snm{Zhang}\ead[label=e1]{zhangxz24@mails.jlu.edu.cn}},
	\author[A,B]{\inits{S.}\fnms{Yong}~\snm{Li}\ead[label=e2]{liyong@jlu.edu.cn}},
	\author[C]{\inits{S.}\fnms{Kaizhi}~\snm{Wang}\ead[label=e3]{kzwang@sjtu.edu.cn}}.
	\address[A]{School of Mathematics, Jilin University, ChangChun, People's Republic of China}
	
	\address[B]{Center for Mathematics and Interdisciplinary Sciences,\\    ~~~ Northeast Normal University, ChangChun, People's Republic of China\\
		}
	\address[C]{School of Mathematical Sciences, CMA-Shanghai, Shanghai Jiao Tong University, Shanghai 200240, China \printead{e1,e2,e3}}
\end{aug}

\begin{abstract}
This paper investigates the application of KAM theory to the stochastic nonlinear Schr{\"{o}}dinger equation on infinite lattices, focusing on the stability of low-dimensional invariant tori in the sense of most probable paths. For generality, we provide an abstract proof within the framework of stochastic Hamiltonian systems on infinite lattices. We begin by constructing the Onsager-Machlup functional for these systems in a weighted infinite sequence space. Using the Euler-Lagrange equation, we identify the most probable transition path of the system's trajectory under stochastic perturbations. Additionally, we establish a large deviation principle for the system and derive a rate function that quantifies the deviation of the system's trajectory from the most probable path, especially in rare events. Combining this with classical KAM theory for the nonlinear Schr{\"{o}}dinger equation, we demonstrate the persistence of low-dimensional invariant tori under small deterministic and stochastic perturbations. Furthermore, we prove that the probability of the system's trajectory deviating from these tori can be described by the derived rate function, providing a new probabilistic framework for understanding the stability of stochastic Hamiltonian systems on infinite lattices.
\end{abstract}



\end{frontmatter}

\section{Introduction}

The Schr{\"{o}}dinger equation, first introduced in 1926 by the renowned Austrian physicist Schr{\"{o}}dinger \cite{61}, serves as a cornerstone of quantum mechanics. It describes the wave-like behavior of particles within a physical system and governs the time evolution of quantum states. In 1973, Hasegawa and Tappert \cite{62,63}, investigating the effects of nonlinear dispersion, derived the nonlinear Schr{\"{o}}dinger equation, which initiated extensive research in this field. As a fundamental nonlinear model that integrates differential equations with matter waves, the nonlinear Schr{\"{o}}dinger equation captures the propagation of wave packets in nonlinear media, including phenomena such as rogue waves in ocean engineering, ultrashort laser pulses, and ion-acoustic waves \cite{64,65,66}. Its applications span various disciplines, including nonlinear optics, radio electronics, biology, telecommunications, and optical soliton communication, making it one of the most actively studied equations in the field of nonlinear differential equations \cite{67,68,69,70,72}.

The classical nonlinear Schrödinger equation admits a Hamiltonian structure \cite{73}. In recent years, significant interest has been generated regarding the application of KAM theory to the nonlinear Schr{\"{o}}dinger equation, with notable contributions by Kuksin \cite{7}, P{\"{o}}schel \cite{8}, Bourgain \cite{9}, Eliasson and Kuksin \cite{10}, Berti and Procesi \cite{75} and others.

Recent studies on stochastic Hamiltonian systems have been flourishing. For example, Wu \cite{17} established a framework for large and moderate deviations, quantifying the probability of rare events in these systems. Talay \cite{18} explored the asymptotic behavior of convergence to an invariant measure. L\'azaro-Cam\'i and Ortega \cite{19} examined the impact of stochastic noise on classical Hamiltonian dynamics. Zhang \cite{20} introduced new computational methods for stochastic flows in Hamiltonian systems based on the Bismut formula. Despite significant progress in this field, a central question remains unresolved: whether quasi-periodic solutions persist in Hamiltonian systems subject to small stochastic perturbations. Recently, in \cite{85}, the first two authors of this paper, Zhang and Li, addressed the case of finite-dimensional stochastic Hamiltonian systems by combining the Onsager-Machlup functional, large deviation theory, and KAM theory. This approach successfully overcame the aforementioned problem, providing the first rigorous proof, in the sense of most probable paths, of the existence of quasi-periodic solutions in finite-dimensional stochastic Hamiltonian systems.

In the present work, we make a substantial extension of this framework: not only do we retain and further refine the analytical strategy developed in \cite{85}, but we also generalize it to infinite lattice Hamiltonian systems, with particular emphasis on the stochastic nonlinear Schr\"{o}dinger equation. Compared with the finite-dimensional case, infinite-dimensional systems exhibit substantially greater analytical complexity and technical obstacles: almost all of the technical lemmas in \cite{85} become inapplicable in this setting, while the analysis of stochastic perturbations on an infinite lattice, together with the derivation of small-ball probability estimates in the infinite-dimensional context, further increases the difficulty of the proof. Confronting these challenges, we develop a new analytical framework for studying the dynamics of infinite-dimensional stochastic Hamiltonian systems in the sense of most probable paths, thereby providing fresh theoretical tools and perspectives for future research in this area.

In this paper, for greater generality, we abstract the system as a general Hamiltonian system on an infinite lattice:
\begin{equation*}
	\begin{cases}
		\,{\rm d} q_i(t) = \frac{\partial H}{\partial p_i}\left( q(t), p(t)\right) \,{\rm d}t, \\
		\,{\rm d} p_i(t) = -\frac{\partial H}{\partial q_i}\left( q(t), p(t)\right) \,{\rm d}t,
	\end{cases} \qquad i \in \mathbb{Z}^m,
\end{equation*}
The corresponding stochastic Hamiltonian system on the infinite lattice can be written in vector form as
\begin{equation}\label{1}
	\begin{cases}
		\,{\rm d} q(t) = \frac{\partial H}{\partial p}\left( q(t), p(t) \right) \,{\rm d}t + \sigma_q(t) \,{\rm d} W_q(t), \\[5pt]
		\,{\rm d} p(t) = -\frac{\partial H}{\partial q}\left( q(t), p(t) \right) \,{\rm d}t + \sigma_p(t) \,{\rm d} W_p(t),
	\end{cases}
\end{equation}
where the infinite-dimensional position and momentum vectors are expanded as  
$
	q(t) = \sum_{i \in \mathbb{Z}^m} q_i(t) \, e_i^{(q)}, 
	\quad 
	p(t) = \sum_{i \in \mathbb{Z}^m} p_i(t) \, e_i^{(p)},
$
where $q_i$ and $p_i$ denote the displacement and momentum, respectively, at lattice site $i$, and $e_i^{(q)}$ and $e_i^{(p)}$ are the canonical unit vectors in the position and momentum directions at site $i$. The Hamiltonian $H(q,p)$ governs the system’s dynamics, with local derivatives $\frac{\partial H}{\partial p_i}$ and $\frac{\partial H}{\partial q_i}$ taken at each site.  
The noise terms $
\sigma_q(t) = \sum_{i \in \mathbb{Z}^m} \sigma_{q_i}(t) \, e_i^{(q)}, 
\sigma_p(t) = \sum_{i \in \mathbb{Z}^m} \sigma_{p_i}(t) \, e_i^{(p)},
$ introduce stochastic perturbations, driven by  
$
W_q(t) = \sum_{i \in \mathbb{Z}^m} W_{q_i}(t) \, e_i^{(q)}, 
W_p(t) = \sum_{i \in \mathbb{Z}^m} W_{p_i}(t) \, e_i^{(p)},
$
where $W_{q_i}(t)$ and $W_{p_i}(t)$ are independent standard Wiener processes defined on a probability space $(\Omega, \mathcal{F}, \mathbb{P})$.

This paper presents a novel analytical framework for investigating the stability of invariant tori in stochastic Hamiltonian systems defined on infinite lattices. The defining conditions are as follows:
\begin{itemize} 
	\item[(C1)] The Hamiltonian function \( H(q,p) \) belongs to \( C^3_b\big(\ell^2_{\rho}(\mathbb{Z}^m; M); \mathbb{R}\big) \), where \( M = \mathbb{T} \times \mathbb{R} \) (The definition of $\ell^2_{\rho}(\mathbb{Z}^m; M)$ is given in Definition \ref{D2.3}). Moreover, the partial derivatives \( \frac{\partial H}{\partial q} \) and \( \frac{\partial H}{\partial p} \) are globally Lipschitz continuous with respect to \( q \) and \( p \); that is, there exists a constant \( L > 0 \) such that, for any \( (q^1,p^1), (q^2,p^2) \in \ell^2_{\rho}(\mathbb{Z}^m; M) \),
	\begin{equation*}
		\left\| \frac{\partial H}{\partial q}(q^1,p^1) -  \frac{\partial H}{\partial q}(q^2,p^2) \right\| 
		\leq L\left( \|q^1 - q^2 \| + \|p^1 - p^2 \| \right),
	\end{equation*}
	\begin{equation*}
		\left\| \frac{\partial H}{\partial p}(q^1,p^1) -  \frac{\partial H}{\partial p}(q^2,p^2) \right\| 
		\leq L\left( \|q^1 - q^2 \| + \|p^1 - p^2 \| \right).
	\end{equation*}
	.
	\item[(C2)]
	The diffusion matrices \( \sigma_p(t) \) and \( \sigma_q(t) \) are block-diagonal matrices, where each diagonal block \( \sigma_{p_i}(t) \) and \( \sigma_{q_i}(t) \) (for each lattice point \( i \in \mathbb{Z}^m \)) is an \( m \times m \) uniformly elliptic matrix that depends continuously on time \( t \). Specifically, they satisfy the following properties:
	\begin{itemize}
		\item[1.] Uniform ellipticity: For all \( t \in [0, T] \), there exists a constant \( \lambda_0 > 0 \) such that for all non-zero vectors \( v \in \mathbb{R}^m \), we have
		\[
		v^\top \sigma_p(t) v \geq \lambda_0 \|v\|^2 \quad \text{and} \quad v^\top \sigma_q(t) v \geq \lambda_0 \|v\|^2.
		\]
		\item[2.] Continuity: The mappings \( t \mapsto \sigma_p(t) \) and \( t \mapsto \sigma_q(t) \) are continuous functions of \( t \) over the interval \( t \in [0, T] \).
	\end{itemize}
\end{itemize}

First, we calculate the Onsager-Machlup functional for the stochastic Hamiltonian system.
\begin{theorem}	\label{T3.1}
	Assume that conditions $(C1)$ and $(C2)$ hold. Let $ (q(t), p(t)) $ be the solution to equation $\eqref{1}$, and let the reference path $ \varphi(t) := \left( \varphi_q(t), \varphi_p(t) \right) $ be a function such that $ \left( (\varphi_q(t), \varphi_p(t)) - (q(0), p(0)) \right) $ belongs to the Cameron-Martin space $ \mathbb{H}^1 $. Then the Onsager-Machlup functional of $ \left( q(t), p(t)\right) $ exists and has the following form:
	\begin{equation*}
		\int_0^T OM(\varphi, \dot{\varphi}) \,{\rm d}t = \int_0^T \left\|\sigma_q^{-1}(s) \left( \dot{\varphi_q}(t) - \frac{\partial H}{\partial \varphi_p}(\varphi_q, \varphi_p) \right)  \right\|^2_{\rho} \,{\rm d}t + \int_0^T \left\| \sigma_p^{-1}(s) \left( \dot{\varphi_p}(t) + \frac{\partial H}{\partial \varphi_p}(\varphi_q, \varphi_p) \right\|^2_{\rho} \right) \,{\rm d}t,
	\end{equation*}
	where $\dot{\varphi}:=\frac{\mathrm{d} \varphi(t)}{\mathrm{d} t} = \left( \frac{\mathrm{d} \varphi_q(t)}{\mathrm{d} t}, \frac{\mathrm{d} \varphi_p(t)}{\mathrm{d} t} \right)$.
\end{theorem}
Since the minimum of the Onsager-Machlup functional $OM(\varphi_q,\varphi_p)$ corresponds to the most probable transition path in the stochastic Hamiltonian system on infinite lattices, we can apply the calculus of variations to Onsager-Machlup functional to obtain this most probable transition path. Exploiting the symplectic structure of the stochastic Hamiltonian system, and working in a weighted space where the position and momentum variables share the same weight, we find that the Onsager-Machlup functional admits a unique representation without correction terms. This allows us to directly identify the most probable transition path $(\varphi_q, \varphi_p)$ as the solution of the deterministic Hamiltonian system
\begin{equation}\label{20}
	\begin{cases}
		d\varphi_q(t) = \frac{\partial H}{\partial \varphi_p}(\varphi_q(t), \varphi_p(t)), \\
		d\varphi_p(t) = -\frac{\partial H}{\partial \varphi_q}(\varphi_q(t), \varphi_p(t)).
	\end{cases}
\end{equation}
To investigate the statistical behavior of the stochastic Hamiltonian system as the noise intensity approaches zero, we rewrite Equation \eqref{1} in the following form for analytical convenience:
\begin{equation}\label{2}
	\begin{cases}
		\mathrm{d}q(t) = \frac{\partial H}{\partial p}(q, p)\,\mathrm{d}t + \epsilon \sigma_q(t)\,\mathrm{d}W_q(t), \\
		\mathrm{d}p(t) = -\frac{\partial H}{\partial q}(q, p)\,\mathrm{d}t + \epsilon \sigma_p(t)\,\mathrm{d}W_p(t),
	\end{cases}
\end{equation}
where \( \epsilon \) represents the noise intensity, with other symbols as defined in \eqref{1}. We then establish a large deviation principle for the most probable path of the stochastic Hamiltonian system.
\begin{theorem} \label{T4.1}
	Assuming that conditions \((C1)\) and \((C2)\) hold, the solution of the stochastic Hamiltonian system \(\eqref{2}\) is given by \( X^\epsilon(t) := (q(t), p(t)) \). As \( \epsilon \to 0 \), the most probable path \( \varphi(t) = (\varphi_q(t), \varphi_p(t)) \) is given by the deterministic Hamiltonian system \(\eqref{20}\). For any path $ X^\epsilon(t) $, the probability that the system deviates from the most probable path satisfies the large deviation principle:
	\begin{equation}
		\epsilon^2 \ln \mathbb{P}(X^\epsilon(t) \in A ) \approx - \inf_{\psi \in A}  J(\psi),
	\end{equation}
	where \( \psi \in A\) denotes an arbitrarily continuous function, and \( A \subset \mathbb{R}^d \) denote an arbitrary measurable set. Furthermore, the rate function \( J(\psi) \) is given by:
	\begin{equation}
		\begin{aligned}
			J(\psi)=
			\begin{cases}
				\frac{1}{2}\left( \int_0^T \left\| \sigma_q^{-1}(t) \left( \dot{\psi}_q - \frac{\partial H}{\partial \psi_p}(\psi_q, \psi_p) \right) \right\|^2_{\rho} \,{\rm d}t\right. \\
				\left. + \int_0^T \left\| \sigma_p^{-1}(t) \left( \dot{\psi}_p + \frac{\partial H}{\partial \psi_q}(\psi_q, \psi_p) \right) \right\|^2_{\rho} \,{\rm d}t \right), & \text{if } \psi - x_0 \in \mathbb{H}^{1};\\
				+\infty, & \text{otherwise}
			\end{cases}
		\end{aligned}
	\end{equation}
	with \( \sigma_q^{-1}(t) \) and \( \sigma_p^{-1}(t) \) being the inverses of the diffusion matrices \( \sigma_q(t) \) and \( \sigma_p(t) \), respectively, and \( \| \cdot \|_\rho \) represents the weighted norm in the space \( L^2 \) with weighting function \( \rho \).
\end{theorem}
These results show that, despite the complexity and uncertainty introduced by stochastic perturbations, the most probable evolution of the system follows the trajectory of the classical Hamiltonian system. Moreover, the determination of the most probable path is independent of the noise intensity $\epsilon$, which only affect the probability of deviations from this path through the rate function in the large deviation principle.

Theorems \ref{T3.1} and \ref{T4.1} are applicable to a broad class of stochastic Hamiltonian systems on arbitrary infinite lattices. In this work, we apply them to the stochastic nonlinear Schr\"odinger equation on an infinite lattice, and, by invoking the classical KAM theorem, establish the persistence of invariant tori along the most probable paths of the stochastic Hamiltonian system.
\begin{theorem}\label{T5.1}
	Consider the stochastic nonlinear Schr{\"{o}}dinger equation
	\begin{equation}\label{50}
		i \frac{\partial u(x,t)}{\partial t} = \frac{\partial^2 u(x,t)}{\partial x^2} - m u(x,t) - f(|u(x,t)|^2) u(x,t) + \epsilon \dot{\eta}(t,x), \quad x \in [0, \pi],~ t \in [0,T].
	\end{equation}
	Assume that the nonlinearity \(f\) is real-analytic and non-degenerate, and that
	\( H := \frac{1}{2} \langle A u, u \rangle + \frac{1}{2} \int_{0}^{\pi} g(|u|^2) \, \mathrm{d}x \)
	satisfies condition \((C1)\), where \(A=-\frac{d^{2}}{dx^{2}}+m\) and \(g(s)=\int_{0}^{s} f(z)\,\mathrm{d}z\).
	Moreover, the noise term \(\eta(t,x)\) can be expanded in the eigenbasis  
	$\phi_j(x) = \sqrt{\frac{2}{\pi}} \sin(jx)$ as  
	$\eta(t, x) = \sum_{j=1}^{\infty} \sigma_j(t) \dot{W}_j(t) \phi_j(x)$,  
	where \(\{W_j(t)\}_{j \geq 1}\) is a family of independent standard one-dimensional Brownian motions, and the coefficients \(\sigma_j(t)\) satisfy assumption \((C2)\).

	Then the most probable transition path of the equation $\eqref{50}$ is governed by the following deterministic equation:
	\[
	i \frac{\partial u(x,t)}{\partial t} = \frac{\partial^2 u(x,t)}{\partial x^2} - m u(x,t) - f(|u(x,t)|^2) u(x,t), \quad x \in [0, \pi],~ t \in [0,T].
	\]
	Moreover, for all \( m \in \mathbb{R} \), all \( n \in \mathbb{N} \), and all \( J = \{j_1 < \cdots < j_n\} \subset \mathbb{N} \), there exists a Cantor manifold \( \mathcal{E}_J \) of real-analytic, linearly stable, Diophantine \( n \)-tori for the equation $\eqref{50}$. This manifold is given by a Lipschitz continuous embedding \( \Phi: T_J[C] \to \mathcal{E}_J \), where the Cantor set \( C \) has full density at the origin. Consequently, \( \mathcal{E}_J \) has a tangent space at the origin equal to \( E_J \), and \( \mathcal{E}_J \) is contained in the space of analytic functions on \( [0, \pi] \).
	
	Furthermore, the probability that the stochastic trajectory \( u^{*}(t,x)\) of equation $\eqref{50}$ deviates from the invariant torus satisfies the following large deviation principle:
	\begin{equation*}
		\epsilon^2 \ln \mathbb{P}(u^{*}(t,x) \in D ) \approx - \inf_{\psi \in D}  J(\psi),
	\end{equation*}
	where \( \psi \in D\) denotes an arbitrarily continuous function, and \( D \subset M \) denote an arbitrary measurable set. Furthermore, the rate function \( J(\psi) \) is given by:
	\begin{equation*}
		\begin{aligned}
			J(\psi)=
			\begin{cases}
				\frac{1}{2} \int_0^T \int_0^{\pi} \left\| \Sigma^{-1}(t) \left(  \frac{\partial \psi(x,t)}{\partial t} + i \left( \frac{\partial^2 \psi(x,t)}{\partial x^2} - m \psi(x,t) - f(|\psi(x,t)|^2) \psi(x,t)\right)  \right) \right\|^2_{a,p}\,{\rm d}x \,{\rm d}t ,& \\&\hspace{-13em} \text{if } \psi - x_0 \in H^1(0,T;\ell^{a,p}(0,\pi));\\
				+\infty, &\hspace{-13em} \text{otherwise},
			\end{cases}
		\end{aligned}
	\end{equation*}
where the inverse of the noise covariance operator $\Sigma(t)$ is defined for a complex-valued function $v(x)$ by
\[
\Sigma^{-1}(t) v := \sum_{j=1}^{\infty} \left( \frac{\langle \operatorname{Re} v, \phi_j \rangle}{\sigma_j^R(t)}  +  i  \frac{\langle \operatorname{Im} v, \phi_j \rangle}{\sigma_j^I(t)} \right) \phi_j(x).
\]
\end{theorem}

This theorem represents a new breakthrough that can be regarded as a stochastic analogue of KAM theory: it establishes the existence of invariant tori in the sense of most probable paths and employs the large deviation principle to quantify the probability of their persistence. In doing so, it extends the classical KAM framework from Hamiltonian systems with deterministic perturbations to those subject to stochastic perturbations.

It is worth emphasizing that the approach of combining the Onsager-Machlup functional, the large deviation principle, and KAM theory--first proposed by the first two authors of this paper--is novel, yet it presents several substantial challenges. For example, achieving an adaptive generalization across different classes of infinite-dimensional Hamiltonian systems remains a major theoretical obstacle. Since many such systems can be discretized into Hamiltonian systems on infinite lattices \cite{9,26}, we begin our analysis within the general framework of Hamiltonian systems on abstract infinite lattices. To resolve convergence issues associated with the Hamiltonian, we introduce a new weighted space $\ell^2_{\rho}$ tailored to the specific requirements of our problem. Furthermore, most analytical tools effective in finite-dimensional settings lose their applicability in the infinite-dimensional context \cite{85}. To overcome these difficulties, we employ a variety of advanced mathematical techniques, including infinite-dimensional Girsanov transforms, Karhunen–Lo\`eve expansions, Skorohod integration, and small-ball probability estimates, thereby providing a refined and robust theoretical foundation for the study of infinite-dimensional Hamiltonian systems.

The structure of this paper is as follows: In Section 2, we review some fundamental definitions of spaces and norms, provide a brief introduction to the Onsager-Machlup functional and large deviation theory, and present several key technical lemmas. In Section 3, we derive the Onsager-Machlup functional for stochastic Hamiltonian systems on infinite lattices and prove that the most probable path corresponds to the solution of the associated deterministic Hamiltonian system. In Section 4, we establish the large deviation principle for stochastic Hamiltonian systems on infinite lattices and obtain the rate function. Finally, in Section 5, we apply the previous results to the one-dimensional stochastic nonlinear Schr{\"{o}}dinger equation on infinite lattices, and, using the corresponding KAM theory, prove the stochastic version of the nonlinear Schr{\"{o}}dinger KAM theory.

\section{Preliminaries}
\subsection{Basic Space}
In this section, we define fundamental spaces and review essential definitions and results related to approximate limits in Wiener spaces (reference \cite{27}).

\begin{definition}
	Let \( \ell^1(\mathbb{Z}^m) \) represent the space of absolutely summable, real-valued sequences indexed by \( \mathbb{Z}^m \), specifically defined as:
	\begin{equation}
		\begin{aligned}
			\ell^1(\mathbb{Z}^m) := \left\lbrace u := \left( u_i\right)_{i \in \mathbb{Z}^m} ~|~ \sum_{i \in \mathbb{Z}^m} |u_i| < \infty \right\rbrace,
		\end{aligned}\notag
	\end{equation}
	where \( |\cdot| \) denotes the absolute value.
\end{definition}

\begin{definition}
	Let \( \ell^2(\mathbb{Z}^m) \) represent the space of square-summable, real-valued sequences indexed by \( \mathbb{Z}^m \), specifically defined as:
	\begin{equation}
		\begin{aligned}
			\ell^2(\mathbb{Z}^m) := \left\lbrace u := \left( u_i\right)_{i \in \mathbb{Z}^m} \in \ell^2 ~|~ \sum_{i \in \mathbb{Z}^m} u_i^{2} < \infty \right\rbrace,
		\end{aligned}\notag
	\end{equation}
	where \( |\cdot| \) denotes the absolute value.
\end{definition}

\begin{definition}
	The \( \ell^\infty(\mathbb{Z}^m) \) space consists of all bounded, real-valued sequences indexed by \( \mathbb{Z}^m \). Specifically, a sequence \( \{b_i\} \) belongs to the \( \ell^\infty(\mathbb{Z}^m) \) space if and only if there exists a constant \( m \) such that \( |b_i| \leq m \) for all \( i \in \mathbb{Z}^m \).
\end{definition}

\begin{definition}
	Let \( \ell^2(\mathbb{Z}^m) \) represent the space of square-summable, real-valued sequences indexed by \( \mathbb{Z}^m \). The standard basis in \( \ell^\infty(\mathbb{Z}^m; M) \), where \( M = \mathbb{T} \times \mathbb{R} \), can be represented by the position basis \( e_i^{(q)} \) and momentum basis \( e_i^{(p)} \) as follows:
	\begin{itemize}
		\item Position basis \( e_i^{(q)} \):
		\begin{equation*}
			e_i^{(q)}(j) = 
			\begin{cases}
				(1, 0), & \text{if } j = i, \\
				(0, 0), & \text{if } j \neq i.
			\end{cases}
		\end{equation*}
		Each \( e_i^{(q)} \) represents the unit vector in the position direction at lattice point \( i \), with a value of \( 1 \) at the position \( i \) and \( 0 \) elsewhere.
		
		\item Momentum basis \( e_i^{(p)} \):
		\begin{equation*}
			e_i^{(p)}(j) = 
			\begin{cases}
				(0, 1), & \text{if } j = i, \\
				(0, 0), & \text{if } j \neq i.
			\end{cases}
		\end{equation*}
		Each \( e_i^{(p)} \) represents the unit vector in the momentum direction at lattice point \( i \), with a value of \( 1 \) at \( i \) in the momentum component and \( 0 \) elsewhere.
	\end{itemize}
	
	Any element \( z = (z_i)_{i \in \mathbb{Z}^m} \in \ell^\infty(\mathbb{Z}^m; M) \) can then be expressed as a linear combination of these basis vectors:
	\[
	z = \sum_{i \in \mathbb{Z}^m} (q_i e_i^{(q)} + p_i e_i^{(p)}),
	\]
	where \( q_i \) and \( p_i \) are the position and momentum components, respectively.
\end{definition}

In this paper, we introduce the weighted space of infinite sequences \( \left( \ell^2_{\rho}(\mathbb{Z}^m; M), \Vert \cdot \Vert_{\rho} \right) \). 

\begin{definition}[Weighted $\ell^2$ metric space on $M=\mathbb T\times\mathbb R$]\label{D2.3-metric}
	Let $\rho=(\rho_i)_{i\in\mathbb Z^m}$ be a sequence of positive weights with $\rho\in \ell^2(\mathbb Z^m)$.
	Define the configuration space
	\[
	\mathcal X_\rho\ :=\ \Bigl\{\,u=(\theta_i,I_i)_{i\in\mathbb Z^m}\ \Big|\ 
	\theta_i\in\mathbb T,\ I_i\in\mathbb R,\ 
	\sum_{i\in\mathbb Z^m}\rho_i^2\,\bigl(d_{\mathbb T}(\theta_i,0)^2+|I_i|^2\bigr)<\infty\Bigr\},
	\]
	where $d_{\mathbb T}$ is the geodesic distance on $\mathbb T$ (taking values in $[0,\pi]$), and $0\in\mathbb T$ is the fixed reference angle.
	
	For $u=(\theta_i,I_i),\,v=(\tilde\theta_i,\tilde I_i)\in\mathcal X_\rho$, define the weighted $\ell^2$ metric
	\[
	d_\rho(u,v)\ :=\ \Biggl(\ \sum_{i\in\mathbb Z^m}\rho_i^2\,\bigl(d_{\mathbb T}(\theta_i,\tilde\theta_i)^2+|I_i-\tilde I_i|^2\bigr)\ \Biggr)^{\!1/2}.
	\]
	Then $(\mathcal X_\rho,d_\rho)$ is a complete separable metric space.
\end{definition}

\begin{definition}\label{D2.3}
	Let the sequence of real numbers \( \left\lbrace \rho_i \right\rbrace_{i \in \mathbb{Z}} \in \ell^2(\mathbb{Z}^m) \cap \ell^1(\mathbb{Z}^m) \). We define a separable Hilbert space \( \left( \ell^2_{\rho}(\mathbb{Z}^m; M), \Vert \cdot \Vert_{\rho} \right) \), where \( M = \mathbb{T} \times \mathbb{R} \), as follows:
	\begin{equation*}
		\ell^2_{\rho}(\mathbb{Z}^m; M) := \left\lbrace u := (q_i, p_i)_{i \in \mathbb{Z}^m} \in M^{\mathbb{Z}^m} ~|~ (q_i, p_i)_{i \in \mathbb{Z}^m} \in \ell^{\infty}, ~ \text{and} ~ \sum_{i \in \mathbb{Z}^m} \left( \rho_i \left\| (q_i, p_i) \right\| \right)^2 < \infty \right\rbrace,
	\end{equation*}
	Where \( \left\| (q_i, p_i) \right\| \) denotes the standard norm in \( M \), defined as \( \left\| (q_i, p_i) \right\| = \sqrt{d_{\mathbb{T}}(q_i, 0)^2 + p_i^2}\), with \( d_{\mathbb{T}} \) representing the geodesic distance on \( \mathbb{T} \).
	
	The norm in this space is defined by:
	\[
	\Vert u \Vert_{\rho} = \left( \sum_{i \in \mathbb{Z}^m} \left( \rho_i \left\| (q_i, p_i) \right\| \right)^2 \right)^{\frac{1}{2}}.
	\]
	For \( u, v \in \ell^2_{\rho}(\mathbb{Z}^m; M) \), the inner product \( \langle u, v \rangle_{\rho} \) is defined as:
	\[
	\langle u, v \rangle_{\rho} = \sum_{i \in \mathbb{Z}^m} \rho_i^2 \left( d_{\mathbb{T}}\left( q_i, \tilde{q}_i\right)  + p_i \tilde{p}_i \right),
	\]
	where \( u = (q_i, p_i)_{i \in \mathbb{Z}^m} \) and \( v = (\tilde{q}_i, \tilde{p}_i)_{i \in \mathbb{Z}^m} \).
\end{definition}

\begin{remark}
	For \( u \in \ell^2_{\rho}(\mathbb{Z}^m) \), we define \( \rho u = \left(\rho_i u_i\right)_{i \in \mathbb{Z}^m} \), and it follows that \( \rho u \in \ell^2(\mathbb{Z}^m) \). Further, we define \( \rho^2 u = \left(\rho_i^2 u_i\right)_{i \in \mathbb{Z}^m} \), and it holds that \( \rho^2 u \in \ell^1(\mathbb{Z}^m) \subseteq \ell^2(\mathbb{Z}^m) \).
\end{remark}

\begin{definition}
	Let \( u(t) \in L^2([0,T], \ell^2_{\rho}(\mathbb{Z}^m)) \). The norm in this space is defined as:
	\begin{displaymath}
		\Vert u\Vert_{L^2_{\rho}} := \Vert u\Vert_{L^2([0,T], \ell^2_{\rho}(\mathbb{Z}^m))} = \left( \int_0^T \Vert u(t) \Vert_{\rho}^2 \, dt \right)^{\frac{1}{2}},
	\end{displaymath}
	and the inner product in \( L^2([0,T], \ell^2_{\rho}(\mathbb{Z}^m)) \) for functions \( u(t) \) and \( v(t) \) is defined by:
	\begin{displaymath}
		\left\langle  u, v \right\rangle_{L^2_{\rho}} := \left\langle  u, v \right\rangle_{L^2([0,T], \ell^2_{\rho}(\mathbb{Z}^m))} = \int_0^T \left\langle  u(t), v(t)\right\rangle_{\rho} \, dt.
	\end{displaymath}
\end{definition}

Let \( W = \left\{ W_t, ~t \in [0, T] \right\} \) be a Brownian motion (Wiener process) defined on the complete filtered probability space \( (\Omega, \mathcal{F}, \left\{ \mathcal{F}_t \right\}_{t \geq 0}, \mathbb{P}) \). Here, \( \Omega \) represents the space of continuous functions vanishing at zero, and \( \mathbb{P} \) denotes the Wiener measure. Let \( \mathbb{H} := L^2([0,T]; \ell^2_{\rho}(\mathbb{Z}^m; M)) \) be a Hilbert space, and let \( \mathbb{H}^1 \) be the Cameron-Martin space defined by:
\begin{displaymath}
	\mathbb{H}^1 := \left\{ f : [0, T] \to \ell^2_{\rho}(\mathbb{Z}^m; M) \in \mathbb{H}^1 ~\big|~f(0) = 0, f ~ \text{is absolutely continuous, and} ~ f^{\prime} \in \mathbb{H} \right\}.
\end{displaymath}
The scalar product in \( \mathbb{H}^1 \) is defined as:
\begin{displaymath}
	\left\langle  f, g \right\rangle_{\mathbb{H}^1} = \left\langle  f^{\prime}, g^{\prime} \right\rangle_{\mathbb{H}} = \left\langle f^{\prime}, g^{\prime} \right\rangle_{L^2_{\rho}}
\end{displaymath}
for all \( f, g \in \mathbb{H}^1 \).

Let \( \mathcal{P}: \mathbb{H}^1 \to \mathbb{H}^1 \) be an orthogonal projection with \( \text{dim}(\mathcal{P} H^1) < \infty \) and the specific expression:
\begin{displaymath}
	\mathcal{P}f = \sum_{i = 1}^{n} \left\langle  h_i, f \right\rangle_{H^1} h_i,
\end{displaymath}
where \( (h_1, ..., h_n) \) is a set of orthonormal basis elements in \( \mathcal{P} \mathbb{H}^1 \). Additionally, we define the \( \mathbb{H}^1 \)-valued stochastic variable:
\begin{displaymath}
	\mathcal{P}^W = \sum_{i = 1}^{n} \bigg( \int_{0}^{T} {h_i^{\prime}} \,{\rm d}W_s \bigg) h_i,
\end{displaymath}
which is independent of the choice of the orthonormal basis $(h_1, \ldots, h_n)$.

\begin{definition}\label{definition 2.1}
	We say that a sequence of orthogonal projections \( \mathcal{P}_n \) on \( \mathbb{H}^1 \) is an approximating sequence of projections if \( \text{dim}(\mathcal{P}_n \mathbb{H}^1) < \infty \) and \( \mathcal{P}_n \) converges strongly to the identity operator \( I \) in \( \mathbb{H}^1 \) as \( n \to \infty \).
\end{definition}

\begin{definition}\label{definition 2.2}
	We say that a semi-norm \( \mathcal{N} \) on \( \mathbb{H}^1 \) is measurable if there exists a stochastic variable \( \tilde{\mathcal{N}} \), satisfying \( \tilde{\mathcal{N}} < \infty \) almost surely, such that for any approximating sequence of projections \( \mathcal{P}_n \) on \( \mathbb{H}^1 \), the sequence \( \mathcal{N}(\mathcal{P}^W_n) \) converges to \( \tilde{\mathcal{N}} \) in probability and \( \mathbb{P}(\tilde{\mathcal{N}} \leq \epsilon) > 0 \) for any \( \epsilon > 0 \). Moreover, if \( \mathcal{N} \) is a norm on \( \mathbb{H}^1 \), then we call it a measurable norm.
\end{definition}

\subsection{Onsager-Machlup Functional}
In the problem of finding the most probable path of a diffusion process, the probability of a single path is zero. Instead, we can search for the probability that the path lies within a certain region, which could be a tube along a differentiable function. This tube is defined as
\[
K(\varphi, \epsilon) = \{ x - x_0 \in \mathbb{H}^1 \mid \varphi - x_0 \in \mathbb{H}^1, \|x - \varphi\| \leq \epsilon, \epsilon > 0 \}.
\]
Once \( \epsilon > 0 \) is given, the probability of the tube can be expressed as
\[
\mu_x(K(\varphi, \epsilon)) = P\left( \{ \omega \in \Omega \mid X_t(\omega) \in K(\varphi, \epsilon) \} \right),
\]
allowing us to compare the probabilities of the tubes for all \( \varphi- x_0 \in \mathbb{H}^1 \), since \( K(\varphi, \epsilon) \in \mathcal{B} \). Here, \( \mu_x \) denotes the probability measure under the initial state \( x \), and \( \mathcal{B} \) is the Borel $\sigma$-field defined on the function space \( \mathbb{H}^1 \), containing all measurable sets generated by open sets.

Thus, the Onsager-Machlup function can be defined as the Lagrangian function that provides the most probable tube. This function plays a crucial role in analyzing path probabilities and rare events in noise-driven systems. Similar to the action functional in classical mechanics, it quantifies the likelihood of various paths within the probabilistic framework.

Onsager and Machlup first introduced this tool in 1953 to describe the probability density of diffusion processes with linear drift and constant diffusion coefficients\cite{28}, \cite{29}. In 1957, Tisza and Manning extended its application to nonlinear equations, while Stratonovich provided a rigorous mathematical framework for the theory in the same year\cite{30}, \cite{31}. In recent years, with deeper research in this area, the Onsager-Machlup functional has been increasingly applied in stochastic systems, particularly for analyzing the most probable path under stochastic perturbations, as discussed in related studies \cite{32,33,34,35,36}.

Based on the above background, we now formally define the Onsager-Machlup function and functional:

\begin{definition}
	Consider a tube surrounding a reference path \( \varphi_t \) with initial value \( \varphi_0 = x \), where \( \varphi_t - x \in \mathbb{H}^1 \). Assuming \( \epsilon \) is small enough, we estimate the probability that the solution process \( X_t \) falls within this tube as:
	\begin{displaymath}
		\mathbb{P} \left\{ \Vert X - \varphi\Vert \leq \epsilon\right\}  \propto C(\epsilon) {\rm exp} \left\{ -\frac{1}{2} \int_{0}^{1} {OM(t, \varphi, \dot{\varphi})} \,{\rm d}t \right\},
	\end{displaymath}
	where \( \propto \) denotes equivalence for small enough \( \epsilon \), and \( \Vert \cdot \Vert \) is an appropriate norm. Here, the integrand \( OM(t, \varphi, \dot{\varphi}) \) is called the Onsager-Machlup function, while the integral \( \int_{0}^{1} {OM(t, \varphi, \dot{\varphi})} \,{\rm d}t \) is called the Onsager-Machlup functional. In the framework of classical mechanics, we also refer to these as the Lagrangian function and the action functional, respectively.
\end{definition}

\subsection{Large deviation principle}
The origins of large deviation theory and its associated research can be traced back to the early 20th century. Cram\'er \cite{37} and Sanov \cite{38} made foundational contributions to the study of large deviations in sequences of independent and identically distributed stochastic variables. Later, Donsker and Varadhan \cite{39} systematically investigated large deviations in the context of Markov processes and explored their relationship with ergodic theory. Their work introduced essential concepts such as Varadhan's integral lemma and the contraction principle, which are not only central results in large deviation theory but also establish profound connections with other areas of mathematics (see \cite{42, 43, 44, 45}). In the 1970s, Freidlin and Wentzell \cite{46} extended this theory to stochastic dynamical systems and stochastic differential equations, particularly in the setting of small perturbations. The Freidlin-Wentzell framework describes the probability of a system deviating from its most likely path and introduces the rate function to quantify the distribution of deviations from typical behavior. Below, we provide the precise definitions of the rate function and the large deviation principle.
\begin{definition}
	A function \( I : E \rightarrow [0, +\infty) \) is called a rate function if \( I \) is lower semicontinuous. Moreover, a rate function \( I \) is called a good rate function if the level set \( \{ x \in E : I(x) \leq K \} \) is compact for each constant \( K < \infty \).
\end{definition}

\begin{definition}
	The stochastic variable sequence \( \{ X^{\epsilon} \} \) is said to satisfy the LDP on \( E \) with rate function \( I \) if the following lower and upper bound conditions hold:
	\begin{itemize}
		\item[(i)] (Lower bound) For any open set \( G \subset E \),
		\[
		\liminf_{\epsilon \rightarrow 0} \epsilon \log \mathbb{P}(X^{\epsilon} \in G) \geq - \inf_{x \in G} I(x).
		\]
		\item[(ii)] (Upper bound) For any closed set \( F \subset E \),
		\[
		\limsup_{\epsilon \rightarrow 0} \epsilon \log \mathbb{P}(X^{\epsilon} \in F) \leq - \inf_{x \in F} I(x).
		\]
	\end{itemize}
\end{definition}

\subsection{KAM Theory}

In Hamiltonian mechanics, invariant tori describe the set of solutions exhibiting quasiperiodic motions. These tori are high-dimensional analogues of closed orbits and arise when a system evolves with incommensurate frequencies. Systems with invariant tori are typically referred to as integrable systems, as their dynamics are regular, confined to these tori, and therefore predictable within the phase space.

The origins of KAM theory lie in the pioneering work of Kolmogorov \cite{81}, Arnold \cite{82}, and Moser \cite{83}, who focused on the stability of these invariant tori under small perturbations. For a nearly integrable Hamiltonian system-where the Hamiltonian is composed of an integrable part plus a small perturbative term-KAM theory asserts that, as long as the perturbation is sufficiently small and certain conditions are met, most of the original invariant tori will persist, albeit with slight deformations. Classic studies on KAM theory can be found in references \cite{47,25,49,50,51,53,57,58}, among others.

We consider small perturbations of an infinite-dimensional Hamiltonian in the parameter-dependent normal form
\[
N = \sum_{1 \leq j \leq n} \omega_j(\xi) y_j + \frac{1}{2} \sum_{j \geq 1} \Omega_j(\xi)(u_j^2 + v_j^2),
\]
defined on the phase space
\[
\mathcal{P}^{a,p} = \mathbb{T}^n \times \mathbb{R}^n \times \ell^{a,p} \times \ell^{a,p} \ni (x, y, u, v),
\]
where \( \mathbb{T}^n \) denotes the usual \( n \)-torus with \( 1 \leq n < \infty \), and \( \ell^{a,p} \) is the Hilbert space of real (later complex) sequences \( w = (w_1, w_2, \dots) \) such that
\[
\|w\|_{a,p}^2 = \sum_{j \geq 1} |w_j|^2 j^{2p} e^{2aj} < \infty,
\]
with \( a \geq 0 \) and \( p \geq 0 \). The frequencies \( \omega = (\omega_1, \dots, \omega_n) \) and \( \Omega = (\Omega_1, \Omega_2, \dots) \) depend on \( n \) parameters \( \xi \in \Pi \subset \mathbb{R}^n \), where \( \Pi \) is a closed, bounded set of positive Lebesgue measure.

To establish the persistence of a large portion of the family of linearly stable rotational tori under small perturbations \( P \) of the Hamiltonian \( N \), as discussed in \cite{25}, we make the following assumptions.

Assumption A: Nondegeneracy. The map $\xi\mapsto\omega(\xi)$ is a lipeomorphism between $\Pi$ and its image, that is, a homeomorphism which is Lipschitz continuous in both directions. Moreover, for all integer vectors $(k,l)\in\mathbb{Z}^n\times\mathbb{Z}^\infty$ with $1\leq|l|\leq2$,
\begin{equation*}
	|\{\xi:\langle k,\omega(\xi)\rangle+\langle l,\Omega(\xi)\rangle = 0\}| = 0
\end{equation*}
and
\begin{equation*}
	\langle l,\Omega(\xi)\rangle\neq0\text{ on }\Pi,
\end{equation*}
where $|\cdot|$ denotes Lebesgue measure for sets, $|l|=\sum_j|l_j|$ for integer vectors, and $\langle\cdot,\cdot\rangle$ is the usual scalar product.

Assumption B: Spectral Asymptotics. There exist $d\geq1$ and $\delta<d - 1$ such that
\begin{equation*}
	\Omega_j(\xi)=j^d+\cdots+O(j^\delta),
\end{equation*}
where the dots stand for fixed lower order terms in $j$, allowing also negative exponents. More precisely, there exists a fixed, parameter-independent sequence $\bar{\Omega}$ with $\bar{\Omega}_j = j^d+\cdots$ such that the tails $\tilde{\Omega}_j=\Omega_j-\bar{\Omega}_j$ give rise to a Lipschitz map
\begin{equation*}
	\tilde{\Omega}:\Pi\to\ell_\infty^{-\delta},
\end{equation*}
where $\ell_\infty^p$ is the space of all real sequences with finite norm $|w|_p=\sup_j|w_j|j^p$.
- Note that the coefficient of $j^d$ can always be normalized to one by rescaling the time. So there is no loss of generality by this assumption. Also, there is no restriction on finite numbers of frequencies.

Assumption C: Regularity. The perturbation $P$ is real analytic in the space coordinates and Lipschitz in the parameters, and for each $\xi\in\Pi$ its hamiltonian vector space field $X_P=(P_y,-P_x,P_v,-P_u)^T$ defines near $\mathcal{T}_0^n$ a real analytic map
\begin{equation*}
	X_P:\mathcal{P}^{a,p}\to\mathcal{P}^{a,\bar{p}},\quad\begin{cases}\bar{p}\geq p&\text{for }d>1,\\\bar{p}>p&\text{for }d = 1.\end{cases}
\end{equation*}
We may also assume that $p-\bar{p}\leq\delta<d - 1$ by increasing $\delta$, if necessary.

To make this quantitative we introduce complex $\mathcal{T}_0^n$-neighbourhoods
\begin{equation*}
	D(s,r):|\text{Im}x|<s,\ |y|<r^2,\ \|u\|_{a,p}+\|v\|_{a,p}<r,
\end{equation*}
where $|\cdot|$ denotes the sup-norm for complex vectors, and weighted phase space norms
\begin{equation*}
	\lVert W\rVert_r=\lVert W\rVert_{\bar{p},r}=|X|+\frac{1}{r^2}|Y|+\frac{1}{r}\|U\|_{a,\bar{p}}+\frac{1}{r}\|V\|_{a,\bar{p}}
\end{equation*}
for $W=(X,Y,U,V)$. Then we assume that $X_P$ is real analytic in $D(s,r)$ for some positive $s$, $r$ uniformly in $\xi$ with finite norm $\lVert X_P\rVert_{r,D(s,r)}=\sup_{D(s,r)}\lVert X_P\rVert_r$, and that the same holds for its Lipschitz semi-norm
\begin{equation*}
	\lVert X_P\rVert_r^{\mathcal{L}}=\sup_{\xi\neq\zeta}\frac{\lVert\Delta_{\xi\zeta}X_P\rVert_r}{|\xi-\zeta|},
\end{equation*}
where $\Delta_{\xi\zeta}X_P=X_P(\cdot,\xi)-X_P(\cdot,\zeta)$, and where the supremum is taken over $\Pi$.

To state the main results we assume that
\begin{equation*}
	|\omega|_{\Pi}^{\mathcal{L}}+|\Omega|_{-\delta,\Pi}^{\mathcal{L}}\leq M<\infty,\quad|\omega^{-1}|_{\omega(\Pi)}^{\mathcal{L}}\leq L<\infty,
\end{equation*}
where the Lipschitz semi-norms are defined analogously to $\lVert X_P\rVert_r^{\mathcal{L}}$. Moreover, we introduce the notations
\begin{equation*}
	\langle l\rangle_d=\max\left(1,\left|\sum j^dl_j\right|\right),\quad A_k = 1 + |k|^\tau,
\end{equation*}
where $\tau\geq n + 1$ is fixed later. Finally, let $\mathcal{Z}=\{(k,l)\neq0,|l|\leq2\}\subset\mathbb{Z}^n\times\mathbb{Z}^\infty$.

\begin{theorem}[\cite{25}]\label{T2.6} 
	Suppose $H = N + P$ satisfies assumptions $A$, $B$ and $C$, and
	\begin{equation*}
		\epsilon=\lVert X_P\rVert_{r,D(s,r)}+\frac{\alpha}{M}\lVert X_P\rVert_{r,D(s,r)}^{\mathcal{L}}\leq\gamma\alpha,
	\end{equation*}
	where $0<\alpha\leq1$ is another parameter, and $\gamma$ depends on $n$, $\tau$ and $s$. Then there exists a Cantor set $\Pi_{\alpha}\subset\Pi$, a Lipschitz continuous family of torus embeddings $\Phi:\mathbb{T}^n\times\Pi_{\alpha}\to\mathcal{P}^{a,p}$, and a Lipschitz continuous map $\omega_*:\Pi_{\alpha}\to\mathbb{R}^n$, such that for each $\xi$ in $\Pi_{\alpha}$, the map $\Phi$ restricted to $\mathbb{T}^n\times\{\xi\}$ is a real analytic embedding of a rotational torus with frequencies $\omega_*(\xi)$ for the hamiltonian $H$ at $\xi$.
	
	Each embedding is real analytic on $|\text{Im}x|<\frac{s}{2}$, and
	\begin{align*}
		\lVert\Phi-\Phi_0\rVert_r+\frac{\alpha}{M}\lVert\Phi-\Phi_0\rVert_r^{\mathcal{L}}&\leq c\epsilon/\alpha,\\
		|\omega_*-\omega|+\frac{\alpha}{M}|\omega_*-\omega|^{\mathcal{L}}&\leq c\epsilon,
	\end{align*}
	uniformly on that domain and $\Pi_{\alpha}$, where $\Phi_0$ is the trivial embedding $\mathbb{T}^n\times\Pi\to\mathcal{T}_0^n$, and $c\leq\gamma^{-1}$ depends on the same parameters as $\gamma$.
	
	Moreover, there exist Lipschitz maps $\omega_{\nu}$ and $\Omega_{\nu}$ on $\Pi$ for $\nu\geq0$ satisfying $\omega_0=\omega$, $\Omega_0=\Omega$ and
	\begin{align*}
		|\omega_{\nu}-\omega|+\frac{\alpha}{M}|\omega_{\nu}-\omega|^{\mathcal{L}}&\leq c\epsilon,\\
		|\Omega_{\nu}-\Omega|_{-\delta}+\frac{\alpha}{M}|\Omega_{\nu}-\Omega|_{-\delta}^{\mathcal{L}}&\leq c\epsilon,
	\end{align*}
	such that $\Pi\backslash\Pi_{\alpha}\subset\bigcup\mathcal{R}_{kl}^{\nu}(\alpha)$, where
	\begin{equation*}
		\mathcal{R}_{kl}^{\nu}(\alpha)=\left\{\xi\in\Pi:\left|\langle k,\omega_{\nu}(\xi)\rangle+\langle l,\Omega_{\nu}(\xi)\rangle\right|<\alpha\frac{\langle l\rangle_d}{A_k}\right\},
	\end{equation*}
	and the union is taken over all $\nu\geq0$ and $(k,l)\in\mathcal{Z}$ such that $|k|>K_02^{\nu-1}$ for $\nu\geq1$ with a constant $K_0\geq1$ depending only on $n$ and $\tau$.
\end{theorem}

\begin{theorem}[\cite{25}]\label{T2.7}
	Let $\omega_{\nu}$ and $\Omega_{\nu}$ for $\nu\geq0$ be Lipschitz maps on $\Pi$ satisfying
	\begin{equation*}
		|\omega_{\nu}-\omega|,\ |\Omega_{\nu}-\Omega|_{-\delta}\leq\alpha,\ |\omega_{\nu}-\omega|^{\mathcal{L}},\ |\Omega_{\nu}-\Omega|_{-\delta}^{\mathcal{L}}\leq\frac{1}{2L},
	\end{equation*}
	and define the sets $\mathcal{R}_{kl}^{\nu}(\alpha)$ as in Theorem A choosing $\tau$ as in (22). Then there exists a finite subset $\mathcal{X}\subset\mathcal{Z}$ and a constant $\tilde{c}$ such that
	\begin{equation*}
		\left|\bigcup_{\substack{(k,l)\notin\mathcal{X}\\\nu}}\mathcal{R}_{kl}^{\nu}(\alpha)\right|\leq\tilde{c}\rho^{n-1}\alpha^{\mu},\quad\mu=\begin{cases}
			1&\text{for }d>1,\\
			\frac{\kappa}{\kappa + 1}&\text{for }d = 1.
		\end{cases}
	\end{equation*}
	for all sufficiently small $\alpha$, where $\rho=\text{diam}\Pi$. The constant $\tilde{c}$ and the index set $\mathcal{X}$ are monotone functions of the domain $\Pi$: they do not increase for closed subsets of $\Pi$. In particular, if $\delta\leq0$, then $\mathcal{X}\subset\{(k,l):0<|k|\leq16LM\}$.
\end{theorem}

The above theorem is based on the framework of KAM theory and applies to certain infinite-dimensional Hamiltonian systems, particularly in the context of weighted infinite lattices or dynamical systems with infinitely many degrees of freedom. The main result is divided into two parts: the analytical part and the geometric part, which are presented as Theorem \ref{T2.6} and Theorem \ref{T2.7}, respectively. The former demonstrates the existence of invariant tori under the assumption that the set of Diophantine frequencies is non-empty. The latter ensures that this condition is indeed satisfied. A detailed discussion of the theorem and its proof can be found in \cite{25}.

\subsection{Karhunen-Lo\`eve expansion}
We calculate here the Karhunen-Lo\`eve expansion for a class of one-dimensional centered mean-square continuous stochastic processes that will appear in the decomposition of $W^Q(t)$.
\begin{definition}
	Assume that stochastic prcess $X : [0,1] \times \Omega \to \mathbb{R}$ is measurable
	for every $t \in [0,1]$. We say stochastic process $X(t, \omega)$ is centered if
	\begin{displaymath}
		\mathbb{E}\left[ X(t, \omega)\right]  = 0 \quad \text{for all}~ t \in [0,1].
	\end{displaymath}
	We say a stochastic process $X(t, \omega)$ is mean-square continuous if
	\begin{displaymath}
		\lim\limits_{\epsilon \to 0} \mathbb{E}\left[ \left( X(t + \epsilon, \omega) - X(t, \omega)\right)^2 \right] = 0 \quad \text{for all}~ t \in [0,1].
	\end{displaymath}
\end{definition}
For a centered mean-square continuous stochastic proces $X(t, \omega)$, we define the integral operator $K : L^2([0,1]) \to L^2([0,1])$ by
\begin{displaymath}
	(Kv)(s) := \int_{0}^{1} k(s,t)v(t)dt, \quad s,t \in [0,1],
\end{displaymath}
where $v(s) \in L^2([0,1])$ and $k(s,t) = \mathbb{E}\left[ X(s, \omega) X(t, \omega)\right] $. So we can show that $K$ is a  compact, positive and self-adjoint operator. According to the spectral theorem, $K$ has a complete set of eigenvectors $\left\lbrace l_i \right\rbrace_{i \in \mathbb{Z}} $ in $L^2([0,1])$ and real non-negative eigenvalues $\left\lbrace \lambda_i^2 \right\rbrace_{i \in \mathbb{Z}} $ (While it is customary to denote the eigenvalues by $\lambda_i$, we use the notation $\lambda_i^2$ here for the sake of convenience.)
:
\begin{displaymath}
	Kl_i = \lambda_i^2 l_i.
\end{displaymath}

Next, we introduce the Karhunen-Lo\`eve expansion theorem related to this paper. For more detailed information on it, please refer to \cite{2}.
\begin{theorem}[\cite{2}]\label{theorem 2.6}
	Let $X : \Omega \times [0,1] \to \mathbb{R}$ be a centered mean-square continuous stochastic process with $X \in L^2(\Omega \times [0,1])$. Then there exists an orthonormal basis $\left\lbrace l_i\right\rbrace_{i \in \mathbb{Z}} $ of $L^2([0,1])$ such that for all $t \in D$,
	\begin{displaymath}
		X(t, \omega) = \sum_{i \in \mathbb{Z}} \lambda_i x_i(\omega) l_i(t),
	\end{displaymath}
	where the coefficients $x_i$ is a sequence of independent, standard normal $\mathcal{N}(0, 1)$ stochastic variables and has the following expression:
	\begin{displaymath}
		x_i(\omega) = \frac{1}{\lambda_i} \int_{0}^{1} X(t, \omega)l_i(t) \mathrm{d} t,
	\end{displaymath}
	and
	\begin{displaymath}
		\lambda_i^2 = \mathrm{Var} \left[ \int_{0}^{1} X(t, \omega)l_i(t) \mathrm{d} t \right].
	\end{displaymath}
\end{theorem}

\subsection{Technical lemmas}
In this section, we will introduce several commonly utilized technical lemmas. Throughout this paper, if not mentioned otherwise, $\mathbb{E} \left(A \big| B\right)$ represents the conditional expectation of $A$ under $B$, and $C$ is a constant and will change with the line.

When we derive the Onsage-Machup functional of SDEs with additive noise, the following lemma is the most basic one, as it ensures that we handle each term separately. Its proof can be found in \cite{54}.
\begin{lemma}[\cite{54}]\label{lemma 2.7}
	For a fixed integer $N \geq 1$, let $X_1, ..., X_N \in \mathbb{R}$ be $N$ stochastic variables defined on $(\Omega, \mathcal{F}, \left\{ \mathcal{F}_t \right\}_{t \geq 0}, \mathbb{P})$ and $\left\{D_{\epsilon}; \epsilon > 0 \right\}$ be a family of sets in $\mathcal{F}$. Suppose that for any $c \in \mathbb{R}$ and any $i = 1, ..., N$, we have
	\begin{displaymath}
		\limsup\limits_{\epsilon \to 0} \mathbb{E}\left({\rm exp}\left\{ c X_i \right\}\big|D_{\epsilon} \right) \leq 1.
	\end{displaymath}
	Then
	\begin{displaymath}
		\limsup\limits_{\epsilon \to 0} \mathbb{E}\left({\rm exp}\left\{ \sum_{i = 1}^{N}c X_i \right\} \big|D_{\epsilon} \right)= 1.
	\end{displaymath}
\end{lemma}
The following two lemmas are about the limit behavior of the expected value of independent, standard normal $\mathcal{N}(0, 1)$ stochastic variables exponential functions, which can be referred to in \cite{55}.
\begin{lemma}[\cite{55}]\label{lmma 2.8}
	Let $(X_i)_{i \in \mathbb{Z}}$ be a sequence of independent, standard normal $\mathcal{N}(0, 1)$ stochastic variables defined on $(\Omega, \mathcal{F}, \left\{ \mathcal{F}_t \right\}_{t \geq 0}, \mathbb{P})$ and let $(\eta_i)_{i \in \mathbb{Z}}$ and $(\theta_i)_{i \in \mathbb{Z}}$ be two real numbers sequences in $ l^2 $. Then
	\begin{displaymath}
		\lim_{\epsilon \to 0} E\left[\exp\left(\sum_{i \in \mathbb{Z}} \eta_i X_i \right) \big| \sum_{i \in \mathbb{Z}} { \theta_i^2 X_i^2} \leq \epsilon\right] = 1.
	\end{displaymath}
	Moreover, for any uniformly bounded stochastic variable $ Y(\omega) $,
	\begin{displaymath}
		\lim_{\epsilon \to 0} E\left[\exp\left(\sum_{i \in \mathbb{Z}} Y(\omega) \eta_i X_i \right) \big| \sum_{i \in \mathbb{Z}} { \theta_i^2 X_i^2} \leq \epsilon\right] = 1.
	\end{displaymath}
\end{lemma}

\begin{lemma}[\cite{55}]\label{lmma 2.9}
	Let $(\eta_i)_{i \in \mathbb{Z}}$ and $(\theta_i)_{i \in \mathbb{Z}}$ be two real numbers sequences in $ l^2 $. And let $(X_i)_{i \in \mathbb{Z}}$ be a sequence of independent, standard normal $\mathcal{N}(0, 1)$ stochastic variables defined on $(\Omega, \mathcal{F}, \left\{ \mathcal{F}_t \right\}_{t \geq 0}, \mathbb{P})$. Assume $ T : l^2 \to l^2$ to be a trace class operator, i.e. $\sum\limits_{i \in \mathbb{Z}} \langle T e_i, e_i \rangle < \infty$ for any orthonormal basis$ \left\lbrace e_i\right\rbrace_{i \in \mathbb{Z}} $ in space $l^2$. Then
	\begin{displaymath}
		\lim_{\varepsilon \to 0} E\left[ \exp\left( \sum_{i,j} X_i X_j T_{ij} \right) \bigg| \sum_{i \in \mathbb{Z}} { \theta_i^2 X_i^2} \leq \epsilon \right] = 1,
	\end{displaymath}
	and more generally
	\begin{displaymath}
		\lim_{\varepsilon \to 0} E\left[ \exp\left( \sum_{i,j} (X_i + \eta_i) X_j T_{ij} \right) \bigg|\sum_{i \in \mathbb{Z}} { \theta_i^2 X_i^2} \leq \epsilon \right] = 1.
	\end{displaymath}
\end{lemma}

The following lemma concerns the small ball probabilities of Gaussian Markov processes under the \( L^p \) norm; for detailed information, refer to \cite{77}. This lemma is formulated for the one-dimensional case, and we will subsequently extend it to an infinite-dimensional version relevant to our study.

Let \( X(t) \), for \( t \in [0,1] \), be a real-valued continuous Gaussian Markov process with mean zero. It is known \cite{78} that the covariance function \( \sigma(s,t) = \mathbb{E}[X(s)X(t)] \) for \( 0 \leq s, t \leq 1 \) satisfies the relation
\[
\sigma(s,t)\sigma(t,u) = \sigma(t,t)\sigma(s,u), \quad 0 \leq s < t < u \leq 1,
\]
which implies the Markov property of \( X(t) \). Consequently, the Gaussian Markov process \( X(t) \) with \( \sigma(s,t) \neq 0 \) for \( 0 < s \leq t < 1 \) can be characterized by
\[
\sigma(s,t) = G(\min(s,t)) H(\max(s,t)),
\]
where \( G(t) > 0 \), \( H(t) > 0 \), and the ratio \( G(t)/H(t) \) is nondecreasing on the interval \( (0,1) \). Moreover, the functions \( G \) and \( H \) are unique up to a constant multiple. We briefly introduce the \( L_p \)-norm on \( C[0,1] \) as follows:
\[
\| f \|_p =
\begin{cases}
	\left( \int_{0}^{1} |f(t)|^p \, dt \right)^{1/p}, & \text{for } 1 \leq p < \infty, \\
	\sup_{0 \leq t \leq 1} |f(t)|, & \text{for } p = \infty.
\end{cases}
\]

\begin{lemma}[\cite{77}]\label{L2.17}
	Let the Gaussian Markov process \(X(t)\) be defined as above. Assume \(H\) and \(G\) are absolutely continuous and \(G/H\) is strictly increasing on the interval \([0,1]\).
	
	If $ \sup_{0<t\leq 1} H(t)<\infty$, or $H(t)$ is nonincreasing in a neighborhood of $0$, then
	\[
	\lim_{\epsilon \to 0} \epsilon^{2} \ln \mathbb{P}(\|X(t)\|_p \leq \epsilon)=-\kappa_p \left(\int_{0}^{1} (G'H - H'G)^{p/(2 + p)} dt\right)^{(2 + p)/p},
	\]
	where
	\[
	\kappa_p = 2^{2/p} p (\lambda_1(p)/(2 + p))^{(2 + p)/p} 
	\]
	and
	\[
	\lambda_1(p)=\inf\left\{\int_{-\infty}^{\infty} \|x\|^p \phi^2(x) dx+\frac{1}{2} \int_{-\infty}^{\infty} (\phi'(x))^2 dx\right\}>0,
	\]
	the infimum is taken over all \(\phi \in L_2(-\infty,\infty)\) such that \(\int_{-\infty}^{\infty} \phi^2(x) dx = 1\).
\end{lemma}
\begin{corollary}[\cite{77}]\label{C2.19}
	Let \( f \) be a locally bounded Borel-measurable function on \([0, \infty)\) such that \( f \in L^2(\mathbb{R}_+) \), and let \( \{B(s)\}_{s \geq 0} \) denote a standard Brownian motion. Then, as established in \cite{84}, the process
	\[
	Z(t) = \int_{0}^{t} f(s) \, \mathrm{d}B(s), \quad t \geq 0,
	\]
	is a Gaussian Markov process with covariance function
	\[
	\sigma_Z(s, t) = \mathrm{Cov}(Z(s), Z(t)) = \int_{0}^{\min(s, t)} f^2(u) \, \mathrm{d}u, \quad s, t \geq 0.
	\]
	Consequently, for \( 1 \leq p \leq \infty \), the following limit holds:
	\[
	\lim_{\epsilon \to 0} \epsilon^2 \ln \mathbb{P} \left( \int_{0}^{\infty} |Z(t)|^p \, \mathrm{d}t \leq \epsilon^p \right) = -\kappa_p \left( \int_{0}^{\infty} |f(t)|^{2p/(2 + p)} \, \mathrm{d}t \right)^{(2 + p)/p},
	\]
	where \( \kappa_p \) is the constant defined in lemma \ref{L2.17}.
\end{corollary}

To adapt to the context of the paper, we define \( X(t) \) as follows:
\[
X(t) := \int_{0}^{t} \sigma(s) \,\mathrm{d} W(s),
\]
where \( \sigma(t) \) is a time-dependent infinite-dimensional diagonal matrix, specifically given by
\[
\sigma(t) = \text{diag} \left( \sigma_1(t), \sigma_2(t), \sigma_3(t), \dots \right)
\]
with each \( \sigma_j(t) \) being a function of time \( t \). Additionally, there exist positive constants \( m \) and \( M \) such that \( m < \sigma_j(t) < M \) for all \( j \in \mathbb{Z} \), $t \in [0,1]$. The process \( W(s) = \left\{ W_j(s), j \in \mathbb{Z}, s \in [0, 1] \right\} \) represents an infinite-dimensional Brownian motion, where each component \( W_j(s) \) is an independent standard one-dimensional Brownian motion.

We extend Corollary $\ref{C2.19}$ to the following form:
\begin{lemma}\label{L2.18}
	For the infinite-dimensional stochastic process \( X(t) \), we have the following results:
	\begin{align}
		\lim_{\epsilon \to 0} \epsilon^{2} \ln \mathbb{P}\left( \left\| X(t) \right\|_{L^2_{\rho}} \leq \epsilon \right) \geq -\kappa_2 C_\rho^2 M^2,
	\end{align}
	where $ \sum_{j \in \mathbb{Z}} \rho_j = C_\rho $ and \( \kappa_2 \) is the constant defined in lemma \ref{L2.17}.
\end{lemma}
\begin{proof}
	Firstly, as established in \cite{84}, the process
	\[
	X(t) = \int_{0}^{t} \sigma(s) \,\mathrm{d} W(s),
	\]
	is a Gaussian Markov process with covariance function
	\[
	\mathrm{Cov}(X(s), X(t)) = \int_{0}^{\min(s, t)}  \sigma^2(u) \, \mathrm{d}u, \quad s, t \geq 0.
	\]
	by the definition of \( L^2_{\rho} \), we have
	\begin{align*}
		\mathbb{P}\left( \left\| X(t) \right\|_{L^2_{\rho}} \leq \epsilon \right) 
		&= \mathbb{P}\left( \left(\int_{0}^{1}  \sum_{j \in \mathbb{Z}}  \left(\rho_j \int_{0}^{t} \sigma_j(s) \,\mathrm{d} W_j(s) \right)^2 \,\mathrm{d}t \right)^{\frac{1}{2}} \leq \epsilon \right)\\
		&= \mathbb{P}\left( \int_{0}^{1}  \sum_{j \in \mathbb{Z}} \rho_j^2 \left( \int_{0}^{t} \sigma_j(s) \,\mathrm{d} W_j(s)\right)^{2} \,\mathrm{d}t \leq \epsilon^2 \right).
	\end{align*}

	Due to the properties of the weights in the infinite-dimensional space we have defined, the sum $ \sum_{j \in \mathbb{Z}} \rho_j$ is finite. Without loss of generality, we set $ \sum_{j \in \mathbb{Z}} \rho_j = C_\rho$. Therefore, we define events A and B as follow
	\begin{align*}
		&A = \left\lbrace \omega \in \Omega | \int_{0}^{1}  \sum_{j \in \mathbb{Z}} \rho_j \left( \int_{0}^{t} \rho_j^{\frac{1}{2}} \sigma_j(s) \,\mathrm{d} W_j(s)\right)^{2} \,\mathrm{d}t \leq \epsilon^2 \right\rbrace,\\
		&B = \left\lbrace \omega \in \Omega | \int_{0}^{1} \left( \int_{0}^{t} \rho_j^{\frac{1}{2}} \sigma_j(s) \,\mathrm{d} W_j(s)\right)^{2} \,\mathrm{d}t \leq \frac{\epsilon^2}{C_\rho} ,\quad \forall j \in \mathbb{Z} \right\rbrace.
	\end{align*}
	For a centered Gaussian process, the variance structure allows us to compare the magnitudes of events $A$ and $B$. Clearly, event $B$ is contained in event $A$, which immediately yields $\mathbb{P}(A) \geq \mathbb{P}(B)$. Moreover, since the components $W_j(t)$ are independent, the probability of event $B$ factorizes as a product by the multiplication rule:
	\begin{align*}
		\mathbb{P}\left( \left\| X(t) \right\|_{L^2_{\rho}} \leq \epsilon \right) 
		&\geq \prod_{j \in \mathbb{Z}} \mathbb{P}\left( \int_{0}^{1}   \left( \int_{0}^{t} \rho_j^{\frac{1}{2}} \sigma_j(s) \,\mathrm{d} W_j(s) \right)^2 \,\mathrm{d}t  \leq \frac{\epsilon^2}{C_\rho} \right).
	\end{align*}
	By taking the logarithm and applying Lemma $\ref{L2.17}$ in this case, we obtain
	\begin{align*}
		\epsilon^{2} \ln \mathbb{P}\left( \left\| X(t) \right\|_{L^2_{\rho}} \leq \epsilon \right) 
		&\geq  \sum_{j \in \mathbb{Z}} \epsilon^{2} \ln \mathbb{P}\left( \int_{0}^{1}   \left( \int_{0}^{t} \rho_j^{\frac{1}{2}}\sigma_j(s) \,\mathrm{d} W_j(s) \right)^2 \,\mathrm{d}t  \leq \frac{\epsilon^2}{C_\rho} \right)\\
		&\geq -\kappa_2 C_\rho \sum_{j \in \mathbb{Z}} \rho_j \left( \int_0^1 | \sigma_j(t)| dt \right)^{2}\\
		&\geq -\kappa_2 C_\rho^2 M^2.
	\end{align*}
\end{proof}

\section{Proof of Theorem \ref{T3.1}}
The main objective of this section is to present the complete proof of Theorem \ref{T3.1}. Our approach relies on deriving the Onsager-Machlup functional for Hamiltonian systems on infinite lattices by computing the ratio of probabilities of perturbed paths within a small neighborhood of a reference path. The key tools used in this derivation include Girsanov’s theorem, the Karhunen–Loève expansion, small-ball probability estimates, and several technical lemmas from Section 2.6. It is worth noting that our results hold for any finite time interval $[0,T]$. For clarity of exposition, however, we set $T=1$ and restrict the discussion to the interval $[0,1]$. The extension to a general interval $[0,T]$ is straightforward and can be obtained through an appropriate time rescaling, which does not affect the validity of the results.

\begin{proof}[proof of Theorem \ref{T3.1}]

	Given a reference path \( \varphi(t) = (\varphi_q(t), \varphi_p(t)) \), which is a deterministic continuous path, and \( (\varphi_q(t), \varphi_p(t)) - (q(0), p(0)) \in \mathbb{H}^1 \), we define the perturbed integral equation with respect to the reference path \( \varphi(t) \) as follows:
	\begin{equation}\label{3}
		\begin{cases}
			\Phi_q(t) = \varphi_q(t) + \int_0^t \sigma_q(s) \,{\rm d}W_q(s),\\
			\Phi_p(t) = \varphi_p(t) + \int_0^t \sigma_p(s) \,{\rm d}W_p(s).
		\end{cases}
	\end{equation}
	To eliminate the drift terms in the original equation, we apply infinite-dimensional Girsanov's theorem (see \cite{56} Theorem 10.14) by introducing a new probability measure \( \tilde{\mathbb{P}} \). Under \( \tilde{\mathbb{P}} \), the transformed Brownian motions are defined as follows:
	\begin{equation}\label{4}
		\begin{aligned}
			\tilde{W}_q(t) &= W_q(t) - \int_0^t \sigma_q^{-1}(s) \left( \frac{\partial H}{\partial \Phi_p}(\Phi_q, \Phi_p) - \dot{\varphi}_q(s) \right) \,{\rm d}s,\\
			\tilde{W}_p(t) &= W_p(t) - \int_0^t \sigma_p^{-1}(s) \left( -\frac{\partial H}{\partial \Phi_q}(\Phi_q, \Phi_p) - \dot{\varphi}_p(s) \right) \,{\rm d}s.
		\end{aligned}
	\end{equation}
	Under the conditions \( (C1) \) and \( (C2) \), we can show that the Novikov condition
	\[
	\mathbb{E}^{\mathbb{P}}\left( \exp\left\{ \frac{1}{2} \int_{0}^{1} \left\| \frac{\partial H}{\partial \Phi_p}(\Phi_q, \Phi_p) - \dot{\varphi}_q(s) \right\|_{\rho}^2 + \left\| \frac{\partial H}{\partial \Phi_q}(\Phi_q, \Phi_p) + \dot{\varphi}_p(s) \right\|_{\rho}^2 \, dt \right\} \right) < + \infty
	\]
	is clearly satisfied. By Girsanov's theorem, it follows that \( \tilde{W}_q(t) \) and \( \tilde{W}_p(t) \) are \( n \)-dimensional standard Brownian motions under the new probability measures \( \tilde{\mathbb{P}}_q \) and \( \tilde{\mathbb{P}}_p \), respectively. Substituting the Brownian motions defined in Equation \(\eqref{4}\) into Equation \(\eqref{3}\), we obtain:
	\[
	\begin{cases}
		\,{\rm d} \Phi_q(t) = \frac{\partial H}{\partial \Phi_p}(\Phi_q, \Phi_p) \,{\rm d}t + \sigma_q(t) \,{\rm d}\tilde{W}_q(t),\\
		\,{\rm d} \Phi_p(t) = -\frac{\partial H}{\partial \Phi_q}(\Phi_q, \Phi_p) \,{\rm d}t + \sigma_p(t) \,{\rm d}\tilde{W}_p(t).
	\end{cases}
	\]
	It follows that, under the new joint probability measure \( \tilde{\mathbb{P}} := \tilde{\mathbb{P}}_q \otimes \tilde{\mathbb{P}}_p \), \( (\Phi_q(t), \Phi_p(t)) \) is a solution to the stochastic Hamiltonian system \(\eqref{1}\).
	
	To facilitate the transformation between two measures, we define the Radon-Nikodym derivative \(\mathcal{R} := \frac{d\tilde{\mathbb{P}}}{d\mathbb{P}}\), which represents the change in measure from \(\mathbb{P}\) to \(\tilde{\mathbb{P}}\). This derivative is given by an exponential martingale associated with the drift term, describing the behavior of Brownian motion under the new measure after eliminating the drift. For the position variable $q$ and the momentum variable $p$, the Radon-Nikodym derivatives are respectively:
	\begin{displaymath}
		\begin{aligned}
			\frac{d\tilde{\mathbb{P}_q}}{d\mathbb{P}_q} &= \exp\left( \int_0^1 \left\langle \sigma_q^{-1}(s) \left( \frac{\partial H}{\partial \Phi_p}(\Phi_q, \Phi_p) - \dot{\varphi}_q(s) \right), \,{\rm d}W_q(s) \right\rangle_{\rho} - \frac{1}{2} \int_0^1 \left\| \sigma_q^{-1}(s) \left( \frac{\partial H}{\partial \Phi_p}(\Phi_q, \Phi_p) - \dot{\varphi}_q(s) \right) \right\|^2_{\rho} \,{\rm d}s \right),\\
			\frac{d\tilde{\mathbb{P}}_p}{d\mathbb{P}_p} &= \exp\left( \int_0^1 \left\langle \sigma_p^{-1}(s) \left( -\frac{\partial H}{\partial \Phi_q}(\Phi_q, \Phi_p) - \dot{\varphi}_p(s) \right), \,{\rm d}W_p(s) \right\rangle_{\rho} - \frac{1}{2} \int_0^1 \left\| \sigma_p^{-1}(s) \left( \frac{\partial H}{\partial \Phi_q}(\Phi_q, \Phi_p) + \dot{\varphi}_p(s) \right) \right\|^2_{\rho} \,{\rm d}s \right).
		\end{aligned}
	\end{displaymath}
	So,
	\begin{displaymath}
		\begin{aligned}
			\mathcal{R} &= \exp\left( \int_0^1 \left\langle \sigma_q^{-1}(s) \left( \frac{\partial H}{\partial \Phi_p}(\Phi_q, \Phi_p) - \dot{\varphi}_q(s) \right), \,{\rm d}W_q(s) \right\rangle_{\rho} - \int_0^1 \left\langle \sigma_p^{-1}(s) \left( \frac{\partial H}{\partial \Phi_q}(\Phi_q, \Phi_p) + \dot{\varphi}_p(s) \right), \,{\rm d}W_p(s) \right\rangle_{\rho} \right. 
			\\ &\quad \left. - \frac{1}{2} \int_0^1 \left\| \sigma_q^{-1}(s) \left( \frac{\partial H}{\partial \Phi_p}(\Phi_q, \Phi_p) - \dot{\varphi}_q(s) \right) \right\|^2_{\rho} \,{\rm d}s - \frac{1}{2} \int_0^1 \left\| \sigma_p^{-1}(s) \left( \frac{\partial H}{\partial \Phi_q}(\Phi_q, \Phi_p) + \dot{\varphi}_p(s) \right) \right\|^2_{\rho} \,{\rm d}s \right).
		\end{aligned}
	\end{displaymath}
	To simplify the notation, let \( W^{\sigma}(t) := \left( W^{\sigma}_{q}(t), W^{\sigma}_{p}(t) \right) \) with
	\begin{displaymath}
		\begin{aligned}
			W^{\sigma}_{q}(t) := \int_0^t \sigma_q(s) \, dW_q(s), \quad W^{\sigma}_{p}(t) := \int_0^t \sigma_q(s) \, dW_q(s).
		\end{aligned}
	\end{displaymath}
	By applying Girsanov's theorem, we can naturally obtain the probability that the path of the stochastic Hamiltonian system \eqref{1} remains close to the reference path \( \varphi(t) \). This probability can be expressed as:
		\begin{align}\label{9}
			&\quad \frac{\mathbb{P}\left(\left\|  (q, p) - (\varphi_q, \varphi_p) \right\|_{L^2_{\rho}}  \leq \epsilon\right)}{\mathbb{P}\left(\left\|  W^{\sigma} \right\|_{L^2_{\rho}}  \leq \epsilon\right)}
			= \frac{\tilde{\mathbb{P}}\left(\left\| (\Phi_q, \Phi_p) -(\varphi_q, \varphi_p) \right\|_{L^2_{\rho}} \leq \epsilon\right)}{\mathbb{P}\left(\left\| W^{\sigma} \right\|_{L^2_{\rho}} \leq \epsilon\right)} 
			= \frac{\mathbb{E} \left( \mathcal{R}\mathbb{I}_{ \left\| W^{\sigma} \right\|_{L^2_{\rho}} \leq \epsilon} \right)}{\mathbb{P}\left( \left\| W^{\sigma} \right\|_{L^2_{\rho}} \leq \epsilon\right)} \notag
			\\&= \mathbb{E}\left( \mathcal{R} \big| \left\| W^{\sigma} \right\|_{L^2_{\rho}} \leq \epsilon \right) 
			\\& = \exp\left\lbrace -\frac{1}{2} \left(  \int_0^1 \left\| \sigma_q^{-1}(t) \left( \dot{\varphi}_q - \frac{\partial H}{\partial \varphi_p}(\varphi_q, \varphi_p) \right) \right\|^2_{\rho} \,{\rm d}t + \int_0^1 \left\| \sigma_p^{-1}(t) \left( \dot{\varphi}_p + \frac{\partial H}{\partial \varphi_q}(\varphi_q, \varphi_p) \right) \right\|^2_{\rho} \,{\rm d}t \right) \right\rbrace \notag
			\\& \quad \times \mathbb{E} \left( \exp\left\lbrace  \sum_{i=1}^{6} B_i \right\rbrace  \bigg| \left\| W^{\sigma} \right\|_{L^2_{\rho}} \leq \epsilon \right), \notag
		\end{align}
	where $B_i$ represents the deviations in the path arising from drift and disturbances, it exhibits stochastic properties. This is elaborated upon in the following detailed expression:
	\begin{align*}
		B_1 &= \int_0^1 \left\langle \sigma_q^{-1}(s) \frac{\partial H}{\partial \Phi_p}(\Phi_q, \Phi_p) , \,{\rm d}W_q(s) \right\rangle_{\rho} - \int_0^1 \left\langle \sigma_p^{-1}(s)  \frac{\partial H}{\partial \Phi_q}(\Phi_q, \Phi_p), \,{\rm d}W_p(s) \right\rangle_{\rho},\\
		B_2 &= - \int_0^1 \left\langle \sigma_q^{-1}(s) \dot{\varphi}_q(s), \,{\rm d}W_q(s) \right\rangle_{\rho} - \int_0^1 \left\langle \sigma_p^{-1}(s) \dot{\varphi}_p(s), \,{\rm d}W_p(s) \right\rangle_{\rho},\\
		B_3 &= \frac{1}{2} \int_0^1 \left\| \sigma_q^{-1}(s) \frac{\partial H}{\partial \varphi_p}(\varphi_q, \varphi_p)  \right\|^2_{\rho} \,{\rm d}s - \frac{1}{2} \int_0^1 \left\| \sigma_q^{-1}(s)  \frac{\partial H}{\partial \Phi_p}(\Phi_q, \Phi_p) \right\|^2_{\rho} \,{\rm d}s,\\
		B_4 &= \frac{1}{2} \int_0^1 \left\| \sigma_p^{-1}(s) \frac{\partial H}{\partial \varphi_q}(\varphi_q, \varphi_p)  \right\|^2_{\rho} \,{\rm d}s - \frac{1}{2} \int_0^1 \left\| \sigma_p^{-1}(s)  \frac{\partial H}{\partial \Phi_q}(\Phi_q, \Phi_p) \right\|^2_{\rho} \,{\rm d}s,\\
		B_5 &= \int_0^1 \left\langle \sigma_q^{-2}(s) \left( \frac{\partial H}{\partial \Phi_p}(\Phi_q, \Phi_p) - \frac{\partial H}{\partial \varphi_p}(\varphi_q, \varphi_p) \right), \dot{\varphi}_q(s) \right\rangle_{\rho} \,{\rm d}s,\\
		B_6 &= - \int_0^1 \left\langle \sigma_p^{-2}(s) \left( \frac{\partial H}{\partial \Phi_q}(\Phi_q, \Phi_p) - \frac{\partial H}{\partial \varphi_q}(\varphi_q, \varphi_p) \right), \dot{\varphi}_p(s) \right\rangle_{\rho} \,{\rm d}s.
	\end{align*}

	First, we know that \( W^{\sigma}(t) = \left( W^{\sigma}_q(t), W^{\sigma}_p(t) \right) \) is a centered mean-square continuous stochastic process. To analyze these processes, we define the covariance operators for the position and momentum components, denoted by \( K_q \) and \( K_p \), respectively. These operators describe the covariance structure of the processes in the weighted space \(L^2([0, 1]; \ell^2_{\rho})\), taking into account the weighting function \(\rho\). Specifically, the covariance operators \( K_q \) and \( K_p \) are defined by their respective covariance functions \( k_q(s, t) \) and \( k_p(s, t) \) as follows:
	\[
	(K_q v)(s) := \int_0^1 k_q(s, t) v(t) \, \mathrm{d}t, \quad (K_p v)(s) := \int_0^1 k_p(s, t) v(t) \, \mathrm{d}t, \quad s, t \in [0,1],
	\]
	where
	\[
	k_q(s, t) = \mathbb{E}[\rho W^{\sigma}_q(s) \cdot \rho W^{\sigma}_q(t)] = \int_0^{\min(s, t)} \rho^2 \sigma_q(u)^2 \, \mathrm{d}u,
	\]
	and
	\[
	k_p(s, t) = \mathbb{E}[\rho W^{\sigma}_p(s) \cdot \rho W^{\sigma}_p(t)] = \int_0^{\min(s, t)} \rho^2 \sigma_p(u)^2 \, \mathrm{d}u.
	\]
	These operators are compact, positive, and self-adjoint, allowing us to utilize the Karhunen-Lo\`eve expansion for our stochastic processes. According to the spectral theorem, each of the covariance operators \( K_q \) and \( K_p \) possesses a complete set of eigenfunctions \(\{l_{j,i}^{(q)}\}_{j \in \mathbb{Z}}\) and \(\{l_{j,i}^{(p)}\}_{j \in \mathbb{Z}}\) in \(L^2([0,1])\) for the position and momentum components, respectively. They are associated with real, non-negative eigenvalues \(\{\lambda_{j,i}^{(q)}\}_{j \in \mathbb{Z}}\) and \(\{\lambda_{j,i}^{(p)}\}_{j \in \mathbb{Z}}\), such that for each component we have
	\[
	K_q l_{j,i}^{(q)} = \lambda_{j,i}^{(q)} l_{j,i}^{(q)} \quad \text{and} \quad K_p l_{j,i}^{(p)} = \lambda_{j,i}^{(p)} l_{j,i}^{(p)}.
	\]
	This spectral decomposition allows the Karhunen-Lo\`eve expansion to be applied separately to the position and momentum components \( q \) and \( p \), facilitating their independent treatment in the analysis. Based on the above content and Theorem \ref{theorem 2.6}, we present the Karhunen-Lo\`eve expansion of \(\rho W^{\sigma}(t)\) as follows:
	\[
	\begin{aligned}
		\rho W^{\sigma}(t) &= \sum_{i \in \mathbb{Z}^m} \left( \left( \int_{0}^{t} {\rho}_{i} \sigma_{q_i}(s) \, dW_{q_i}(s) \right) e^{(q)}_i + \left( \int_{0}^{t} {\rho}_{i} \sigma_{p_i}(s) \, dW_{p_i}(s) \right) e^{(p)}_i \right)\\
		&= \sum_{i \in \mathbb{Z}^m} \sum_{j \in \mathbb{Z}} \lambda_{j,i}^{(q)} x_{j,i}^{(q)} \left( l_{j,i}^{(q)}(t) \otimes e_i^{(q)} \right)(t) +  \sum_{i \in \mathbb{Z}^m} \sum_{j \in \mathbb{Z}} \lambda_{j,i}^{(p)} x_{j,i}^{(p)} \left( l_{j,i}^{(p)}(t) \otimes e_i^{(p)} \right)(t),
	\end{aligned}
	\]
	where \( \left\lbrace l_{j,i}^{(q)}(t) \right\rbrace_{j \in \mathbb{Z}}\) and \( \left\lbrace l_{j,i}^{(p)}(t) \right\rbrace_{j \in \mathbb{Z}}\) are orthogonal bases of \( L^2([0,1]) \), and \(\lambda_{j,i}^{(q)}\), \( x_{j,i}^{(q)} \), and \( l_{j,i}^{(q)}(t) \); \(\lambda_{j,i}^{(p)}\), \( x_{j,i}^{(p)} \), and \( l_{j,i}^{(p)}(t) \) are the eigenvalues, independent stochastic variables, and eigenfunctions corresponding to the position component \( q \) and momentum component \( p \), respectively. The specific expressions are as follows:
	\[
	\begin{aligned}
		& \lambda_{j,i}^{(q)} = \left( \mathrm{Var} \left[ \int_0^1 \left( \int_0^t \rho_{i} \sigma_{q_i}(s) \, \mathrm{d}W_{q_i}(s) \right) l_{j,i}^{(q)}(t) \, \mathrm{d}t \right]\right)^{\frac{1}{2}},\\
		&\lambda_{j,i}^{(p)} = \left( \mathrm{Var} \left[ \int_0^1 \left( \int_0^t \rho_{i} \sigma_{p_i}(s) \, \mathrm{d}W_{p_i}(s) \right) l_{j,i}^{(p)}(t) \, \mathrm{d}t \right]\right)^{\frac{1}{2}},\\
		&x_{j,i}^{(q)}(\omega) = \frac{1}{\lambda_{j,i}^{(q)}} \int_0^1 \left( \int_0^t \rho_{i} \sigma_{q_i}(s) \, \mathrm{d}W_{q_i}(s) \right) l_{j,i}^{(q)}(t) \, \mathrm{d}t,\\
		&x_{j,i}^{(p)}(\omega) = \frac{1}{\lambda_{j,i}^{(p)}} \int_0^1 \left( \int_0^t \rho_{i} \sigma_{p_i}(s) \, \mathrm{d}W_{p_i}(s) \right) l_{j,i}^{(p)}(t) \, \mathrm{d}t.
	\end{aligned}
	\]
	Notably, \(\{l_{j,i}^{(q)} \otimes e_i^{(q)}, l_{j,i}^{(p)} \otimes e_i^{(p)}; {j \in \mathbb{Z}}, i \in \mathbb{Z}^m\}\) forms an orthonormal basis for \( L^2([0, 1]; \ell^{\infty}) \), ensuring that for any \({j \in \mathbb{Z}}\) and \( i \in \mathbb{Z}^m \),
	\[
	\text{Cov}((W^{\sigma_q}, l_{j,i}^{(q)} \otimes e_i^{(q)})_{L^2([0,1]; \ell^{\infty})}) = (\lambda_{j,i}^{(q)})^2, \quad \text{and} \quad \text{Cov}((W^{\sigma_p}, l_{j,i}^{(p)} \otimes e_i^{(p)})_{L^2([0,1]; \ell^{\infty})}) = (\lambda_{j,i}^{(p)})^2.
	\]
	Thus, we can express the norm \(\Vert W^Q \Vert^2_{L^2_{\rho}}\) as
	\[
	\Vert W^Q \Vert^2_{L^2_{\rho}} = \sum_{i \in \mathbb{Z}^m} \sum_{j \in \mathbb{Z}} \left((\lambda_{j,i}^{(q)})^2 (x_{j,i}^{(q)})^2 + (\lambda_{j,i}^{(p)})^2 (x_{j,i}^{(p)})^2\right).
	\]
	Since \(\{\sigma_{q_i}(t)\}_{i \in \mathbb{Z}^m}\), \(\{\sigma_{p_i}(t)\}_{i \in \mathbb{Z}^m} \in L^2([0, 1]; \ell^2_{\rho})\), we conclude that \(\sum\limits_{i \in \mathbb{Z}^m} \sum\limits_{j \in \mathbb{Z}} \left((\lambda_{j,i}^{(q)})^2 + (\lambda_{j,i}^{(p)})^2\right) < +\infty\).

	For the second term $B_2$, we have
	\begin{equation*}
		\begin{aligned}
			B_2 &= - \int_0^1 \left\langle \sigma_q^{-1}(s) \dot{\varphi}_q(s), \,{\rm d}W_q(s) \right\rangle_{\rho} - \int_0^1 \left\langle \sigma_p^{-1}(s) \dot{\varphi}_p(s), \,{\rm d}W_p(s) \right\rangle_{\rho}
			\\&= - \sum_{i \in \mathbb{Z}^m} \left( \int_0^1 \rho_{i}^2 \sigma_{q_i}^{-1}(s) \dot{\varphi}_{q_i}(s) \,{\rm d}W_{q_i}(s) + \int_0^1 \rho_{i}^2 \sigma_{p_i}^{-1}(s) \dot{\varphi}_{p_i}(s) \,{\rm d}W_{p_i}(s) \right).
		\end{aligned}
	\end{equation*}
	It is straightforward to demonstrate that the sequences \( \{\sigma_{q,i}^{-1}(s) \dot{\varphi}_{q,i}(s)\}_{i \in \mathbb{Z}^m} \in L^2([0,1]; l^2_{\rho}) \) and \( \{\sigma_{p,i}^{-1}(s) \dot{\varphi}_{p,i}(s)\}_{i \in \mathbb{Z}^m} \in L^2([0,1]; l^2_{\rho}) \). Since for any \( i \in \mathbb{Z}^m \), the sets \( \left\lbrace l_{j,i}^{(q)}(t) \right\rbrace_{j \in \mathbb{Z}} \) and \( \left\lbrace l_{j,i}^{(p)}(t) \right\rbrace_{j \in \mathbb{Z}} \) form orthogonal bases in \( L^2([0,1]) \), we can project \(\rho_{i}^2 \sigma_{q_i}^{-1}(s) \dot{\varphi}_{q_i}(s)\) and \(\rho_{i}^2 \sigma_{p_i}^{-1}(s) \dot{\varphi}_{p_i}(s)\) onto these bases:
		\begin{align*}
			B_2 &= - \sum_{i \in \mathbb{Z}^m} \left( \int_0^1 \sum_{j \in \mathbb{Z}}  \pi_{j,i}^{(q)} l_{j,i}^{(q)}(s) \,{\rm d}W_{q_i}(s) + \int_0^1 \sum_{j \in \mathbb{Z}}  \pi_{j,i}^{(p)} l_{j,i}^{(p)}(s) \,{\rm d}W_{p_i}(s) \right)\\
			&= - \sum_{i \in \mathbb{Z}^m} \sum_{j \in \mathbb{Z}} \left( \int_0^1  \pi_{j,i}^{(q)} l_{j,i}^{(q)}(s) \,{\rm d}W_{q_i}(s) + \int_0^1   \pi_{j,i}^{(p)} l_{j,i}^{(p)}(s) \,{\rm d}W_{p_i}(s) \right)\\
			&= - \sum_{i \in \mathbb{Z}^m} \sum_{j \in \mathbb{Z}} \left(  \pi_{j,i}^{(q)} I^q_i \left( l_{j,i}^{(q)}\right) + \pi_{j,i}^{(p)} I^p_i \left( l_{j,i}^{(p)}\right) \right),
		\end{align*}
	where \( \pi_{j,i}^{(q)} \) and \( \pi_{j,i}^{(p)} \) represent the projection coefficients corresponding to the position component \( q \) and the momentum component \( p \), respectively. The specific expressions are as follows:
	\[
	\pi_{j,i}^{(q)} = \left( \mathrm{Var} \left[ \int_0^1 \rho_{i}^2  \sigma_{q_i}^{-1}(s) \dot{\varphi}_{q_i}(s) l_{j,i}^{(q)}(s) \, \mathrm{d}s \right] \right)^{\frac{1}{2}}, \quad 
	\pi_{j,i}^{(p)} = \left( \mathrm{Var} \left[ \int_0^1 \rho_{i}^2  \sigma_{p_i}^{-1}(s) \dot{\varphi}_{p_i}(s) l_{j,i}^{(p)}(s) \, \mathrm{d}s \right] \right)^{\frac{1}{2}}.
	\]
	In addition, \( I^q_i \left( l_{j,i}^{(q)}(s) \right) := \int_0^1  l_{j,i}^{(q)}(s) \,{\rm d}W_{q_i}(s) \) and \( I^p_i \left( l_{j,i}^{(p)}(s) \right) := \int_0^1  l_{j,i}^{(p)}(s) \,{\rm d}W_{p_i}(s) \) are two sequences of independent, standard normal stochastic variables \( \mathcal{N}(0, 1) \). Consequently, the result follows directly from Lemma \ref{lmma 2.8}:
	\begin{equation}\label{10}
		\limsup\limits_{\epsilon \to 0} \mathbb{E}\left( {\rm exp} \left\{ c B_2 \right\} \big| \left\| W^{\sigma} \right\|_{L^2_{\rho}} < \epsilon \right) = 1
	\end{equation}
	for all \( c \in \mathbb{R} \). 
	
	For the third term $B_3$, 
	\begin{displaymath}
		\begin{aligned}
			B_3 &= \frac{1}{2} \int_0^1 \left\| \sigma_q^{-1}(s) \frac{\partial H}{\partial \varphi_p}(\varphi_q, \varphi_p)  \right\|^2_{\rho} \,{\rm d}s - \frac{1}{2} \int_0^1 \left\| \sigma_q^{-1}(s)  \frac{\partial H}{\partial \Phi_p}(\Phi_q, \Phi_p) \right\|^2_{\rho} \,{\rm d}s
			\\ &\leq \frac{1}{2} \int_{0}^{1} {\sigma_q^{-2}(s) \left\|  \frac{\partial H}{\partial \varphi_p}(\varphi_q, \varphi_p) - \frac{\partial H}{\partial \Phi_p}(\Phi_q, \Phi_p) \right\|^2_{\rho} }\\
			&+ {2 \sigma_q^{-2}(s) \left\| \frac{\partial H}{\partial \varphi_p}(\varphi_q, \varphi_p) - \frac{\partial H}{\partial \Phi_p}(\Phi_q, \Phi_p) \right\|_{\rho} \left\| \frac{\partial H}{\partial \Phi_p}(\Phi_q, \Phi_p) \right\|_{\rho}} \,{\rm d}s
			\\ &\leq \frac{1}{2} \int_{0}^{1} {\sigma_q^{-2}(s) \left\|  \frac{\partial H}{\partial \varphi_p}(\varphi_q, \varphi_p) - \frac{\partial H}{\partial \Phi_p}(\Phi_q, \Phi_p) \right\|^2_{\rho}} \,{\rm d}s\\
			&+ \int_{0}^{1} {\sigma_q^{-2}(s) \left\|\frac{\partial H}{\partial \varphi_p}(\varphi_q, \varphi_p) - \frac{\partial H}{\partial \Phi_p}(\Phi_q, \Phi_p) \right\|_{\rho} \left\| \frac{\partial H}{\partial \Phi_p}(\Phi_q, \Phi_p) \right\|_{\rho}} \,{\rm d}s.
		\end{aligned}
	\end{displaymath}
	In Condition $(C1)$, since \(\frac{\partial H}{\partial p}\) is Lipschitz continuous, we have the following estimate:  
	\begin{equation}\label{11}
		\left\| \frac{\partial H}{\partial \Phi_p}(\Phi_q, \Phi_p)  -  \frac{\partial H}{\partial \varphi_p}(\varphi_q, \varphi_p)\right\|_{\rho} = \left\| \frac{\partial H}{\partial \left( \varphi_p + W^{\sigma}_p \right)}((\varphi_q + W^{\sigma}_q), (\varphi_p + W^{\sigma}_p))-  \frac{\partial H}{\partial \varphi_p}(\varphi_q, \varphi_p)\right\|_{\rho}  \leq L \left\| W^{\sigma} \right\|_{\rho}.
	\end{equation}
	Inequality $\eqref{11}$, H{\"{o}}lder's inequality and the boundedness of $\frac{\partial H}{\partial \Phi_p}(\Phi_q, \Phi_p)$ and $\sigma_q^{-1}(t)$ imply that
	\begin{displaymath}
		\begin{aligned}
			B_3 \leq C \int_{0}^{1} \left\| W^{\sigma} \right\|_{\rho} \,{\rm d}s
			\leq C \left( \int_{0}^{1} \left\| W^{\sigma} \right\|_{\rho}^2 \,{\rm d}s\right)^{\frac{1}{2}} \left( \int_{0}^{1} 1 \,{\rm d}s\right)^{\frac{1}{2}}
			\leq C \left\| W^{\sigma} \right\|_{L^2_{\rho}}.
		\end{aligned}
	\end{displaymath}
	Thus,
	\begin{equation}\label{12}
		\limsup\limits_{\epsilon \to 0} \mathbb{E}\left({\rm exp}\left\{ cB_3 \right\} \big|\left\| W^{\sigma} \right\|_{L^2_{\rho}} < \epsilon \right) = 1
	\end{equation}
	for all $c \in \mathbb{R}$.
	
	For the fourth term \( B_4 \), under Condition \( (C1) \), since \( \frac{\partial H}{\partial q} \) is Lipschitz continuous, we obtain an inequality similar to inequality \(\eqref{11}\):
	\begin{equation}
		\left\| \frac{\partial H}{\partial \Phi_q}(\Phi_q, \Phi_p) - \frac{\partial H}{\partial \varphi_q}(\varphi_q, \varphi_p) \right\|_{\rho} \leq L \left\| W^{\sigma} \right\|_{\rho}.\notag
	\end{equation}
	Using an approach analogous to that applied for the third term \( B_3 \), we obtain
	\begin{equation}\label{13}
		\limsup_{\epsilon \to 0} \mathbb{E}\left( \exp\left\{ cB_4 \right\} \,\big|\, \left\| W^{\sigma} \right\|_{L^2_{\rho}} < \epsilon \right) = 1
	\end{equation}
	for all \( c \in \mathbb{R} \).
	
	For the fifth term \( B_5 \), applying inequality $\eqref{11}$, H{\"{o}}lder's inequality and the boundedness of $\dot{\varphi}_q(t)$ and $\sigma_q^{-1}(t)$, we have
	\begin{displaymath}
		\begin{aligned}
			B_5 &= \int_0^1 \left\langle \sigma_q^{-2}(s) \left( \frac{\partial H}{\partial \Phi_p}(\Phi_q, \Phi_p) - \frac{\partial H}{\partial \varphi_p}(\varphi_q, \varphi_p) \right), \dot{\varphi}_q(s) \right\rangle_{\rho} \,{\rm d}s\\
			& \leq C \int_0^1 \left\|  \frac{\partial H}{\partial \Phi_p}(\Phi_q, \Phi_p) - \frac{\partial H}{\partial \varphi_p}(\varphi_q, \varphi_p) \right\|_{\rho} \,{\rm d}s\\
			& \leq C\int_0^1 L \left\| W^{\sigma} \right\|_{\rho} \,{\rm d}s\\
			& \leq C\left\| W^{\sigma} \right\|_{L^2_{\rho}}.
		\end{aligned}
	\end{displaymath}
	Thus,
	\begin{equation}\label{14}
		\limsup\limits_{\epsilon \to 0} \mathbb{E}\left({\rm exp}\left\{ cB_5 \right\} \big|\left\|  W^{\sigma} \right\|_{L^2_{\rho}} < \epsilon \right) = 1
	\end{equation}
	for all $c \in \mathbb{R}$.
	
	For the sixth term $B_6$, using an analogous approach to the one applied for the fifth term $B_5$, we obtain
	\begin{equation}\label{15}
		\limsup\limits_{\epsilon \to 0} \mathbb{E}\left({\rm exp}\left\{ cB_6 \right\} \big|\left\|  W^{\sigma} \right\|_{L^2_{\rho}} < \epsilon \right) = 1
	\end{equation}
	for all $c \in \mathbb{R}$.
	
	For the first term $ B_1 $, to combine its two parts into a cohesive expression and simplify the notation, we define the following vectorized terms. Let the inverse noise matrix be:
	\[
	\sigma^{-1}(t) := \begin{bmatrix}\sigma_q^{-1}(t) & 0 \\ 0 & \sigma_p^{-1}(t) \end{bmatrix},
	\]
	where \(\sigma_q^{-1}(t) =  \sum\limits_{i \in \mathbb{Z}^m} \sigma^{-1}_{q_i}(t) e_i^{(q)} \) and \(\sigma_p^{-1}(t) = \sum\limits_{i \in \mathbb{Z}^m} \sigma^{-1}_{p_i}(t) e_i^{(p)} \) represent the inverse noise components associated with the position and momentum, respectively. We also define the system's gradient vector as:
	\[
	H'(y) := \begin{pmatrix}\frac{\partial H}{\partial \Phi_p}(\Phi_q, \Phi_p) \\ - \frac{\partial H}{\partial \Phi_q}(\Phi_q, \Phi_p) \end{pmatrix},
	\]
	where \( \Phi = (\Phi_q, \Phi_p) \) is the state vector of the system. Additionally, let \( W(t) \) and the differential increment \(\,{\rm d}W(t)\) be given by:
	\[
	W(t) := \begin{pmatrix} W_q(t) \\ W_p(t) \end{pmatrix}, \quad \,{\rm d}W(t) := \begin{pmatrix} \,{\rm d}W_q(t) \\ \,{\rm d}W_p(t) \end{pmatrix}.
	\]
	Under the assumption of small perturbations, it is feasible to apply a Taylor series expansion to $H'(\Phi)$. Specifically, we have
	\begin{displaymath}
		\begin{aligned}
			H'(\Phi) &= \begin{pmatrix}\frac{\partial H}{\partial \varphi_p}(\varphi_q, \varphi_p) \\ - \frac{\partial H}{\partial \varphi_q}(\varphi_q, \varphi_p)\\ \end{pmatrix} + \begin{bmatrix}\frac{\partial^2H}{\partial \varphi_q \partial \varphi_p}(\varphi_q, \varphi_p) & \frac{\partial^2H}{\partial^2 \varphi_p}(\varphi_q, \varphi_p) \\ - \frac{\partial^2H}{\partial^2 \varphi_q}(\varphi_q, \varphi_p) & - \frac{\partial^2H}{\partial \varphi_p \partial \varphi_q}(\varphi_q, \varphi_p) \\ \end{bmatrix} \begin{pmatrix} W^{\sigma}_q(t) \\ W^{\sigma}_p(t)  \\ \end{pmatrix}
			+  \begin{pmatrix}R_q(t) \\ R_p(t)  \\ \end{pmatrix}
			\\& := H'(\varphi) + J_{H'}(\varphi) W^{\sigma} + R(t).
		\end{aligned}
	\end{displaymath}
	The term \( J_{H'}(\varphi) \) is referred to as the infinite-dimensional Jacobian matrix or the Fréchet derivative of \( H' \) at \( \varphi \), while \( R(t) \) represents the higher-order remainder term. Based on the properties of the Taylor series expansion, when \( H \in C^3_b \) and \( \left\| W^{\sigma} \right\|_{L^2_{\rho}} \leq \epsilon \), we can estimate the remainder term \( R(t) \) as follows:
	\begin{displaymath}
		\left\|  R(t) \right\|_{L^2_{\rho}}  \leq k \epsilon^2,
	\end{displaymath}
	where \( k \) is a finite positive constant.
	
	Hence, $B_1$ can be written as:
	\begin{displaymath}
		\begin{aligned}
			B_1 &= \int_0^1 \left\langle \sigma_q^{-1}(s) \frac{\partial H}{\partial \Phi_p}(\Phi_q, \Phi_p) , \,{\rm d}W_q(s) \right\rangle_{\rho} - \int_0^1 \left\langle \sigma_p^{-1}(s)  \frac{\partial H}{\partial \Phi_q}(\Phi_q, \Phi_p), \,{\rm d}W_p(s) \right\rangle_{\rho}\\
			&= \int_0^1 \left\langle \sigma^{-1}(s)  H'(\Phi),  \,{\rm d}W(s)\right\rangle_{\rho} \\
			&= \int_0^1 \left\langle \sigma^{-1}(s)  H'(\varphi),  \,{\rm d}W(s)
			\right\rangle_{\rho} + \int_0^1 \left\langle \sigma^{-1}(s)  J_{H'}(\varphi) W^{\sigma},  \,{\rm d}W(s)\right\rangle_{\rho} 
			+ \int_0^1 \left\langle \sigma^{-1}(s)   R(s),  \,{\rm d}W(s)\right\rangle_{\rho} 
			\\ &:= B_{11} + B_{12} + B_{13}.
		\end{aligned}
	\end{displaymath}

	The term \( B_{11} \) has a form similar to that of \( B_2 \):
	\begin{displaymath}
		\begin{aligned}
			B_{11} =  \int_0^1 \left\langle \sigma^{-1}(s)  H'(\varphi),  \,{\rm d}W(s)\right\rangle_{\rho}.
		\end{aligned}
	\end{displaymath}
	Due to $H \in C^3_b$, we can show that $\sigma^{-1}(s)  H'(\varphi) \in L^2([0, 1]; \ell^2_{\rho})$. Using the same method as item $B_2$ yields
	\begin{equation}\label{16}
		\limsup\limits_{\epsilon \to 0} \mathbb{E}\left({\rm exp}\left\{ cB_{11} \right\} \big|\left\|  W^{\sigma} \right\|_{L^2_{\rho}}  < \epsilon \right) = 1
	\end{equation}
	for all $c \in \mathbb{R}$.
	
	The term $B_{12}$ is a double stochastic integral with respect to $W$:
	\begin{equation*}
		\begin{aligned}
			B_{12} &= \int_0^1 \left\langle \sigma^{-1}(s)  J_{H'}(\varphi) W^{\sigma},  \,{\rm d}W(s)\right\rangle_{\rho}\\
			&= \int_0^1 \left\langle \rho \sigma^{-1}(s)  J_{H'}(\varphi) \rho W^{\sigma},  \,{\rm d}W(s)\right\rangle\\
			&= \sum_{k \in \mathbb{Z}^m} \int_0^1 \rho_{k} \sigma_{q_k}^{-1}(s) \left\langle J_{H'}(\varphi) \sum_{i \in \mathbb{Z}^m} \sum_{j \in \mathbb{Z}} \lambda_{j,i}^{(q)} x_{j,i}^{(q)} l_{j,i}^{(q)}(s) \otimes e_i^{(q)}, e^{(q)}_k \right\rangle \,{\rm d}W_{q_k}(s)\\
			&\quad + \sum_{k \in \mathbb{Z}^m} \int_0^1 \rho_{k} \sigma_{p_k}^{-1}(s) \left\langle J_{H'}(\varphi)  \sum_{i \in \mathbb{Z}^m} \sum_{j \in \mathbb{Z}} \lambda_{j,i}^{(p)} x_{j,i}^{(p)} l_{j,i}^{(p)}(s) \otimes e_i^{(p)}, e^{(p)}_k \right\rangle \,{\rm d}W_{p_k}(s)\\
			&= \sum_{k \in \mathbb{Z}^m} \sum_{i \in \mathbb{Z}^m} \sum_{j \in \mathbb{Z}} \int_0^1 \rho_{k} \sigma_{q_k}^{-1}(s) \lambda_{j,i}^{(q)} x_{j,i}^{(q)}  \left\langle J_{H'}(\varphi) l_{j,i}^{(q)}(s) \otimes e_i^{(q)}, e^{(q)}_k \right\rangle \,{\rm d}W_{q_k}(s)\\
			&\quad + \sum_{k \in \mathbb{Z}^m} \sum_{i \in \mathbb{Z}^m} \sum_{j \in \mathbb{Z}} \int_0^1 \rho_{k} \sigma_{p_k}^{-1}(s) \lambda_{j,i}^{(p)} x_{j,i}^{(p)} \left\langle J_{H'}(\varphi)  l_{j,i}^{(p)}(s) \otimes e_i^{(p)}, e^{(p)}_k \right\rangle \,{\rm d}W_{p_k}(s)\\
			& := B_{12}^{(q)} + B_{12}^{(p)}.
		\end{aligned}
	\end{equation*}
	Here, we consider \( B_{12}^{(q)} \). Let
	\begin{displaymath}
		h^{(q)}_{j,i}(s) := \frac{1}{\lambda_{j,i}} \int_{s}^{1} \sigma_{q_i}(t) l^{(q)}_{j,i}(t) \,{\rm d}t
	\end{displaymath}
	for $i \in \mathbb{Z}^m$ and $j \in \mathbb{Z}$. Then, for any $i \in \mathbb{Z}^m$, the set $\left\lbrace h^{(q)}_{j,i}(s)\right\rbrace_{j \in \mathbb{Z}} $ forms an orthonormal basis of $L^2([0, 1])$. Notice that
	\begin{displaymath}
		\begin{aligned}
			\langle h^{(q)}_{j_1,i}, h^{(q)}_{j_2,i} \rangle
			&= \frac{1}{\lambda_{j_1,i} \lambda_{j_2,i}} \int_{0}^{1} \left( \int_{s}^{1} \sigma_{q_i}(t) l^{(q)}_{j_1,i}(t) \,{\rm d}t_1\right) \left( \int_{s}^{1} \sigma_{q_i}(t) l^{(q)}_{j_1,i}(t) \,{\rm d}t_2\right) \,{\rm d}s \\ 
			&= \frac{1}{\lambda_{j_1,i} \lambda_{j_2,i}} \int_{0}^{1} \int_{0}^{1} K_i(t, s) l^{(q)}_{j_1,i}(t) l^{(q)}_{j_2,i}(s) \,\mathrm{d}t \,{\rm d}s \\
			&= \langle l^{(q)}_{j_1,i}, l^{(q)}_{j_2,i} \rangle
		\end{aligned}
	\end{displaymath}
	for any $j_1, j_2 \in \mathbb{Z}$, where $K_i(t,s)$ denotes the covariance function defined by
	\begin{displaymath}
		K_i(t, s) = \int_{0}^{t \wedge s} \left( \sigma_{q_i}(u)\right)^2 \,{\rm d}u.
	\end{displaymath}
	To prove that $\left\lbrace h^{(q)}_{j,i} \right\rbrace_{j \in \mathbb{Z}}$ forms a complete basis, one needs to show that if $f(t) \in L^2([0, 1])$ and $\langle h^{(q)}_{j,i}, f \rangle = 0$ for all $j  \in \mathbb{Z}$, then $f(t)$ must be zero. This directly follows from the fact that if 
	
	\begin{displaymath}
		0 = \langle h^{(q)}_{j,i}, f \rangle = \frac{1}{\lambda_{j,i}} \int_{0}^{1} \left( \int_{s}^{1} \sigma_{q_i}(t) l^{(q)}_{j,i}(t) \,{\rm d}t \right) f(s) \,{\rm d}s = \frac{1}{\lambda_{n,i}} \left\langle  l^{(q)}_{j,i}(t), \sigma_{q_i}(t) \int_{0}^{t}  f(s) \,\mathrm{d}s \right\rangle 
	\end{displaymath}
	for all $j  \in \mathbb{Z}$. Since \( l^{(q)}_{j,i}(t) \) is an orthonormal basis in \( L^2[0,1] \), it follows that \( \sigma_{q_i}(t) \int_{0}^{t}  f(s) \,\mathrm{d}s = 0 \), which allows us to obtain $f \equiv 0$. Furthermore,
	\begin{displaymath}
		\begin{aligned}
			x^{(q)}_{j,i} &= \frac{1}{\lambda_{j,i}} \int_{0}^{1} \left( \int_{0}^{t} \sigma_{q_i}(t) \,{\rm d}W_{q_i}(s) \right) l^{(q)}_{j,i}(t) \,{\rm d}t\\
			&= \int_{0}^{1} \left( \frac{1}{\lambda_{j,i}} \int_{s}^{1} \sigma_{q_i}(t) l^{(q)}_{j,i}(t) \,{\rm d}t \right)  \,{\rm d}W_{q_i}(s)\\
			&= I_i(h_{j,i}),
		\end{aligned}
	\end{displaymath}
	where $I^{(q)}_i(h^{(q)}_{j,i}) = \int_{0}^{1}h^{(q)}_{j,i}(s) \,{\rm d}W_{q_i}(s)$.
	
	Let $P$ and $Q$ be two linear operators defined on $L^2([0, 1]; \mathbb{H})$ such that for any $f \in L^2([0, 1])$,
	\begin{displaymath}
		P\left( f(s) \otimes e_i(x) \right):= \rho \sigma(s)^{-1} J_{H'}(\varphi)\left( f(s) \otimes e_i(x) \right),
	\end{displaymath} 
	and $Q(f(s) \otimes e_i(x)) = \left( Q_i f\right) (s) \otimes e_i(x)$ with
	\begin{displaymath}
		(Q_i f)(s) := \int_{s}^{1} \sigma_{q_i}(t) f(t) \,{\rm d}t.
	\end{displaymath}
	So we have
	\begin{displaymath}
		\begin{aligned}
			B_{12}^{(q)} &= \sum_{k \in \mathbb{Z}^m} \sum_{i \in \mathbb{Z}^m} \sum_{j \in \mathbb{Z}} \int_0^1 \rho_{k} \sigma_{q_k}^{-1}(s) \lambda_{j,i}^{(q)} x_{j,i}^{(q)}  \left\langle J_{H'}(\varphi) l_{j,i}^{(q)}(s) \otimes e_i^{(q)}, e^{(q)}_k \right\rangle \,{\rm d}W_{q_k}(s)\\
			&= \sum_{k \in \mathbb{Z}^m} \sum_{i \in \mathbb{Z}^m} \sum_{j \in \mathbb{Z}}  \lambda_{j,i}^{(q)} \int_0^1  x_{j,i}^{(q)}  \left\langle P\left(  l_{j,i}^{(q)}(s) \otimes e_i^{(q)}\right) , e^{(q)}_k \right\rangle \,{\rm d}W_{q_k}(s).
		\end{aligned}
	\end{displaymath}

	It is noteworthy that the stochastic variable $ x_{m,i}$ is measurable with respect to $\mathcal{F}_1$. Hence, the introduction of Skorohod integral is required to handle the anticipating stochastic integrals that arise in the above expression. We transition from It{\^{o}} integrals to Skorohod integrals, capitalizing on the property that they concur on the set $ L^2_a $ of square-integrable adapted processes. For a detailed discussion on Skorohod integrals, refer to \cite{33}. When $k = i$,
	\begin{displaymath}
		\begin{aligned}
			&\int_0^1  x_{j,i}^{(q)}  \left\langle P\left(  l_{j,i}^{(q)}(s) \otimes e_i^{(q)}\right) , e^{(q)}_k \right\rangle \,{\rm d}W_{q_k}(s)\\
			& = \int_0^1  x_{j,i}^{(q)}  \left\langle P\left(  l_{j,i}^{(q)}(s) \otimes e_i^{(q)}\right) , e^{(q)}_i \right\rangle \,{\rm d}W_{q_i}(s)\\
			& = x_{j,i}^{(q)} \sum_{m=1}^{\infty} \left\langle P \left( l_{j,i}^{(q)}(s) \otimes e^{(q)}_i \right) , h^{(q)}_{m,i} \otimes e^{(q)}_i \right\rangle I_i(h^{(q)}_{m,i}) - \left\langle P \left( l_{j,i}^{(q)}(s) \otimes e^{(q)}_i \right) , h^{(q)}_{j,i} \otimes e^{(q)}_i \right\rangle,
		\end{aligned}
	\end{displaymath}
	and when $k \neq i$,
	\begin{displaymath}
		\int_{0}^{1} { x_{j,i}^{(q)} \left\langle P \left( l_{j,i}^{(q)}(s) \otimes e^{(q)}_i \right) , e^{(q)}_k \right\rangle } \,{\rm d}W_{q_k}(t) = x_{j,i}^{(q)} \sum_{m=1}^{\infty} \left\langle P \left( l_{j,i}^{(q)}(s) \otimes e^{(q)}_i \right) , h^{(q)}_{m,k} \otimes e^{(q)}_k \right\rangle I_k(h^{(q)}_{m,k}).
	\end{displaymath}
	Using the fact that $h^{(q)}_{m,k} = \frac{1}{\lambda_{m,k}} Q_k l^{(q)}_{m,k}$ and $I^{(q)}_k(h^{(q)}_{m,k}) = x^{(q)}_{m,k}$, we can write $B^{(q)}_{12}$ in the following way:
	\begin{displaymath}
		\begin{aligned}
			B^{(q)}_{12} &= \sum_{(m,k) \neq (j,i)} \frac{\lambda^{(q)}_{j,i}}{\lambda^{(q)}_{m,k}} x^{(q)}_{m,k} x^{(q)}_{j,i} \left\langle  P (l^{(q)}_{j,i} \otimes e^{(q)}_i), Q \left( l^{(q)}_{m,k} \otimes e^{(q)}_k \right) \right\rangle\\
			&\quad + \sum_{j,i} \left( \left( x^{(q)}_{j,i}\right) ^2 - 1\right)  \left\langle  P (l^{(q)}_{j,i} \otimes e^{(q)}_i), Q \left( l^{(q)}_{j,i} \otimes e^{(q)}_i \right) \right\rangle.
		\end{aligned}
	\end{displaymath}
	Define now the operator $T : l^2_{\mathbb{Z} \times \mathbb{Z}^m} \rightarrow l^2_{\mathbb{Z} \times \mathbb{Z}^m}$ by
	\begin{displaymath}
		T_{(j,i),(m,k)} = \frac{\lambda^{(q)}_{j,i}}{\lambda^{(q)}_{m,k}} \left\langle  P (l^{(q)}_{j,i} \otimes e^{(q)}_i), Q \left( l^{(q)}_{m,k} \otimes e^{(q)}_k \right) \right\rangle, \quad (j,i),(m,k) \in \mathbb{Z} \times \mathbb{Z}^m.
	\end{displaymath}
	Due to
		\begin{align*}
			Q^*Q\left(  l^{(q)}_{m,k} \otimes e_k\right) &= \sum_{n=1}^{\infty} \left\langle (Q^* Q)_k l^{(q)}_{m,k}, l^{(q)}_{n,k} \right\rangle (l^{(q)}_{n,k} \otimes e^{(q)}_k)\\
			&= \sum_{n=1}^{\infty} \left\langle Q_k l^{(q)}_{m,k}, Q_k l^{(q)}_{n,k} \right\rangle (l^{(q)}_{n,k} \otimes e^{(q)}_k) \\
			& = \left( \lambda^{(q)}_{m,k}\right) ^2 \left( l^{(q)}_{m,k} \otimes e^{(q)}_k\right),
		\end{align*}
	and
	\begin{displaymath}
		\begin{aligned}
			\frac{1}{\lambda^{(q)}_{m,k}}Q\left( l^{(q)}_{m,k} \otimes e^{(q)}_k\right) = (h^{(q)}_{m,k} \otimes e^{(q)}_k),
		\end{aligned}
	\end{displaymath}
	we have
	\begin{displaymath}
		\begin{aligned}
			T_{(j,i),(m,k)} &=  \left\langle \lambda^{(q)}_{j,i} P \left( l^{(q)}_{j,i} \otimes e^{(q)}_i\right), \frac{1}{\lambda^{(q)}_{m,k}} Q \left( l^{(q)}_{m,k} \otimes e^{(q)}_k \right) \right\rangle\\
			& =  \left\langle \frac{1}{\lambda^{(q)}_{j,i}} PQ^*Q \left( l^{(q)}_{j,i} \otimes e^{(q)}_i\right) , \frac{1}{\lambda^{(q)}_{m,k}} Q \left( l^{(q)}_{m,k} \otimes e^{(q)}_k \right) \right\rangle\\
			& =  \left\langle PQ^* \left( h^{(q)}_{j,i} \otimes e^{(q)}_i\right) , \left( h^{(q)}_{m,k} \otimes e^{(q)}_k \right) \right\rangle.
		\end{aligned}
	\end{displaymath}
	Then
	\begin{displaymath}
		\begin{aligned}
			B^{(q)}_{12} &= \sum_{(m,k) \neq (j,i)} T_{(j,i),(m,k)} x_{j,i} x_{m,k} + \sum_{j,i} T_{(j,i),(j,i)} (x_{j,i}^2 - 1)\\
			&= \sum_{m,k,j,i} T_{(j,i),(m,k)} x_{j,i} x_{m,k} + \sum_{j,i} T_{(j,i),(j,i)} (x_{j,i}^2 - 2).
		\end{aligned}
	\end{displaymath}
	According to Lemma \ref{lmma 2.9}, we define the self adjoint operator $\tilde{T}=\frac{1}{2} \left(P Q^* + \left( P Q^*\right)^*  \right) $. Due to the boundedness of $B_{12}$, we know that $\tilde{T}$ is a trace class operator. Consequently,
	\begin{displaymath}
		\limsup\limits_{\epsilon \to 0} \mathbb{E}\left(\text{exp} \left\lbrace c\sum_{m,k,j,i} T_{(j,i),(m,k)} x_{j,i} x_{m,k} \right\rbrace  \big|\Vert W^g \Vert_{L^2_{\rho}} < \epsilon \right) = 1
	\end{displaymath}
	for all $c \in \mathbb{R}$. Therefore,
	\begin{equation}
		\begin{aligned}
			&\quad \limsup\limits_{\epsilon \to 0} \mathbb{E}\left({\rm exp}\left\{ B_{12} \right\} \big|\Vert W^g \Vert_{L^2_{\rho}} < \epsilon \right) \\
			& = \limsup\limits_{\epsilon \to 0} \mathbb{E}\left({\rm exp}\left\{ \sum_{j,i} T_{(j,i),(j,i)} (x_{j,i}^2 - 2) \right\} \big|\Vert W^g \Vert_{L^2_{\rho}} < \epsilon \right)\\
			& = {\rm exp}\left\{ - { Tr\left({\tilde{T}} \right)}  \right\}.
		\end{aligned}\notag
	\end{equation}
	Give that $Q(f \otimes e^{(q)}_i) = Q_i(f) \otimes e^{(q)}_i$ with
	\begin{displaymath}
		(Q_i f)(s) := \int_{s}^{1} \sigma_{q_i}(t) f(t) \,{\rm d}t.
	\end{displaymath}
	We have $Q^*(f \otimes e^{(q)}_i) = Q^*_i(f) \otimes e^{(q)}_i$ with
	\begin{displaymath}
		(Q_i^* f)(s) = \sigma_{q_i}(t) \int_{0}^{t}  f(s) \,{\rm d}s.
	\end{displaymath}
	Then
	\begin{displaymath}
		\begin{aligned}
			(PQ^*)(f \otimes e^{(q)}_i)(s) &= \sum_{n=1}^{\infty} \left( \left( \int_0^s f(t) \,{\rm d}t\right) \rho_{n} \sigma_{q_i}(s) \sigma^{-1}_{q_n}(s) \left\langle J_{H'} (\varphi) e^{(q)}_i, e^{(q)}_n  \right\rangle(s) \right)  e^{(q)}_n,
			\\(PQ^*)^*(f \otimes e^{(q)}_i)(s) &= \sum_{n=1}^{\infty} \left( \int_s^1 f(t) \rho_{n} \sigma_{q_i}(t) \sigma^{-1}_{q_n}(t) \left\langle J_{H'} (\varphi) e^{(q)}_i, e^{(q)}_n \right\rangle (t) \,{\rm d}t \right) e^{(q)}_n.
		\end{aligned}\notag
	\end{displaymath}
	We thus have
	\begin{align*}
		&\text{Tr} \left( \frac{1}{2} \left( P Q^* + (P Q^*)^* \right) \right) \\
		&= \sum_{j,i=1}^{\infty} \left\langle \frac{1}{2} \left( P Q^* + (P Q^*)^* \right) \left( {l}^{(q)}_{j,i} \otimes e^{(q)}_i\right), \left( {l}^{(q)}_{j,i} \otimes e^{(q)}_i \right) \right\rangle \\
		&= \frac{1}{2} \sum_{j,i=1}^{\infty} \left\langle  \sum_{n=1}^{\infty} \left( \left( \int_0^s l^{(q)}_{j,i}(t) \, \mathrm{d}t \right) \rho_{n} \sigma_{q_i}(s) \sigma^{-1}_{q_n}(s) \left\langle J_{H'}(\varphi) e^{(q)}_i, e^{(q)}_n \right\rangle \right.\right. \\
		&\quad \left.\left. + \int_s^1 l^{(q)}_{j,i}(t) \rho_{n} \sigma_{q_i}(t) \sigma^{-1}_{q_n}(t) \left\langle J_{H'}(\varphi) e^{(q)}_i, e^{(q)}_n \right\rangle \, \mathrm{d}t \right) e^{(q)}_n, \left( {l}^{(q)}_{j,i} \otimes e^{(q)}_i \right)\right\rangle \\
		&= \frac{1}{2} \sum_{j,i=1}^{\infty} \left\langle \left( \int_0^s l^{(q)}_{j,i}(t) \, \mathrm{d}t \left\langle J_{H'}(\varphi) e^{(q)}_i, e^{(q)}_i \right\rangle_{\rho}   + \int_s^1 l^{(q)}_{j,i}(t) \left\langle J_{H'}(\varphi) e^{(q)}_i, e^{(q)}_i \right\rangle_{\rho}  \, \mathrm{d}t \right),  {l}^{(q)}_{j,i} \right\rangle\\
		&= \frac{1}{2}  \sum_{j,i=1}^{\infty} \left\langle \int_0^1 l^{(q)}_{j,i}(t) \left( \left\langle J_{H'}(\varphi(s)) e^{(q)}_i, e^{(q)}_i \right\rangle_{\rho} \cdot 1_{[0, s]}(t) + \left\langle J_{H'}(\varphi(t)) e^{(q)}_i, e^{(q)}_i \right\rangle_{\rho} \cdot 1_{[s, 1]}(t) \right) \, \mathrm{d}t,  {l}^{(q)}_{j,i} \right\rangle\\
		&= \frac{1}{2}  \sum_{j,i=1}^{\infty} \left\langle \int_0^1 l^{(q)}_{j,i}(t) \hat{K}(s,t) \, \mathrm{d}t,  {l}^{(q)}_{j,i} \right\rangle\\
		&= \frac{1}{2}  \sum_{j,i=1}^{\infty} \int_0^1 \hat{K}(s,s) \, \mathrm{d}s\\
		&= \frac{1}{2} \sum_{i=1}^{\infty}  \int_0^1 \left\langle J_{H'}(\varphi) e^{(q)}_i, e^{(q)}_i \right\rangle_{\rho} \, \mathrm{d}s,
	\end{align*}
	where 
	\begin{equation*}
		\hat{K}(s, t) = \left\langle J_{H'}(\varphi(s \vee t)) e^{(q)}_i, e^{(q)}_i \right\rangle_{\rho}.
	\end{equation*}
	By performing the same derivation for \( B_{12}^{(p)} \) as we did for \( B_{12}^{(q)} \), we obtain the following result: 
	\begin{equation}
		\begin{aligned}
			&\limsup\limits_{\epsilon \to 0} \mathbb{E}\left({\rm exp}\left\{ B^{(q)}_{12} \right\} \big|\Vert W^{\sigma} \Vert_{L^2_{\rho}} < \epsilon \right)
			= {\rm exp}\left\{ - \frac{1}{2} \sum_{i=1}^{\infty}  \int_0^1 \left\langle J_{H'}(\varphi) e^{(q)}_i, e^{(q)}_i \right\rangle_{\rho} \, \mathrm{d}s  \right\},\\
			&\limsup\limits_{\epsilon \to 0} \mathbb{E}\left({\rm exp}\left\{ B^{(p)}_{12} \right\} \big|\Vert W^{\sigma} \Vert_{L^2_{\rho}} < \epsilon \right)
			= {\rm exp}\left\{ - \frac{1}{2} \sum_{i=1}^{\infty}  \int_0^1 \left\langle J_{H'}(\varphi) e^{(p)}_i, e^{(p)}_i \right\rangle_{\rho} \, \mathrm{d}s  \right\}.
		\end{aligned}\notag
	\end{equation}
	Therefore,
	\begin{equation*}
		\begin{aligned}
			\limsup\limits_{\epsilon \to 0} \mathbb{E}\left({\rm exp}\left\{ B_{12} \right\} \big|\Vert W^g \Vert_{L^2} < \epsilon \right) 
			={\rm exp}\left\{ - \frac{1}{2} \sum_{i=1}^{\infty}  \int_0^1 \left\langle J_{H'}(\varphi) e_i, e_i \right\rangle_{\rho} \, \mathrm{d}s  \right\}
			= {\rm exp}\left\{ -\frac{1}{2} \int_{0}^{T} \text{Tr}_{\rho}(  J_{H'}(\varphi) ) \,{\rm d}t  \right\}.
		\end{aligned}
	\end{equation*}
	Given that the Hamiltonian system possesses a symplectic structure and that we apply the same weighting function to both the position and momentum vectors, we can demonstrate that
	\begin{equation*}
		\text{Tr}_{\rho}(  J_{H'}(\varphi) = \text{Tr}_{\rho} \left( \frac{\partial^2H}{\partial \varphi_q \partial \varphi_p}(\varphi_q, \varphi_p) \right)  + \text{Tr}_{\rho} \left( - \frac{\partial^2H}{\partial \varphi_p \partial \varphi_q}(\varphi_q, \varphi_p) \right)  = 0. 
	\end{equation*}
	Thus, we have:
	\begin{equation}\label{17}
		\begin{aligned}
			\limsup\limits_{\epsilon \to 0} \mathbb{E}\left({\rm exp}\left\{ cB_{12} \right\} \big|\Vert W^{\sigma} \Vert_{L^2_{\rho}} < \epsilon \right) = 1
		\end{aligned}
	\end{equation}
	for all $c \in \mathbb{R}$.
	
	Finally, we study the behaviour of the term $B_{13}$. For any $c \in \mathbb{R}$ and $\delta > 0$, we have
	\begin{equation*}
		\begin{aligned}
			\mathbb{E} \left( {\rm exp} \left\{ cB_{13} \right\} \Vert W^{\sigma} \Vert_{L^2_{\rho}} \leq \epsilon \right)
			&=  \int_{0}^{\infty} {e^x \mathbb{P}\left( c\int_0^1 \left\langle \sigma^{-1}(s)   R(s),  \,{\rm d}W(s)\right\rangle_{\rho} > x \big| \Vert W^{\sigma} \Vert_{L^2_{\rho}} \leq \epsilon \right)} \,{\rm d}x\\
			&\leq \int_{\delta}^{\infty} {e^x \mathbb{P}\left( c\int_0^1 \left\langle \sigma^{-1}(s)   R(s),  \,{\rm d}W(s)\right\rangle_{\rho} > x \big| \Vert W^{\sigma} \Vert_{L^2_{\rho}} \leq \epsilon \right)} \,{\rm d}x  + e^{\delta}.
		\end{aligned}
	\end{equation*}
	Define the martingale $M_t = c\int_0^1 \left\langle \sigma^{-1}(s) R(s),  \,{\rm d}W(s)\right\rangle_{\rho} $. We have estimate about its quadratic variation
	\begin{displaymath}
		\langle M_t \rangle = c^2\int_{0}^{t} {\Vert \sigma^{-1}(s) R(s) \Vert^2_{\rho}} \,{\rm d}s \leq C \epsilon^4
	\end{displaymath}
	for some $C > 0$. Using the exponential inequality for martingales, we have
	\begin{displaymath}
		\mathbb{P}\left( c\int_0^1 \left\langle \sigma^{-1}(s)   R(s),  \,{\rm d}W(s)\right\rangle_{\rho} > x, \Vert W^{\sigma} \Vert_{L^2_{\rho}} \leq \epsilon \right) \leq {\rm exp}\left\{ -\frac{x^2}{2C\epsilon^4} \right\}.
	\end{displaymath}
	Based on Lemma \ref{L2.18}, we obtain the following estimate:
	\begin{displaymath}
		\mathbb{P}\left( \left\| W^{\sigma} \right\|_{L^2_{\rho}} \leq \epsilon \right) \geq \exp\left\lbrace - \frac{\kappa_p C_\rho M^2}{\epsilon^2} \right\rbrace.
	\end{displaymath}
	Therefore,
	\begin{equation}\label{18}
		\begin{aligned}
			\mathbb{P}\left(  c\int_0^1 \left\langle \sigma^{-1}(s)   R(s),  \,{\rm d}W(s)\right\rangle_{\rho}  > x \big| \Vert W^{\sigma} \Vert_{L^2_{\rho}} \leq \epsilon \right) &= \frac{\mathbb{P}\left(  c\int_0^1 \left\langle \sigma^{-1}(s)   R(s),  \,{\rm d}W(s)\right\rangle_{\rho}  > x , \Vert W^{\sigma} \Vert_{L^2_{\rho}} \leq \epsilon \right)}{\mathbb{P}\left( \Vert W^{\sigma} \Vert_{L^2_{\rho}} \leq \epsilon \right)} \\
			&\leq C {\rm exp}\left\{ -\frac{x^2}{2C\epsilon^4} +  \frac{ \kappa_2 C_\rho^2 M^2}{\epsilon^2} \right\}.
		\end{aligned}
	\end{equation}
	According to condition $(C3)$, by taking the limit in Inequality \(\eqref{18}\), we obtain
	\begin{equation}\label{19}
		\limsup\limits_{\epsilon \to 0} \mathbb{E}\left({\rm exp}\left\{ cB_{13} \right\} \big|\Vert W^g \Vert_{L^2_\rho} < \epsilon \right) = 1
	\end{equation}
	for all $c \in \mathbb{R}$, as $\epsilon \to 0$ and $\delta \to 0$.
	
	In summary, by Lemma $\ref{lemma 2.7}$ and Inequalities  $\eqref{9}$, $\eqref{10}$, $\eqref{12}$-$\eqref{17}$ and $\eqref{19}$, we have
	\begin{displaymath}
		\begin{aligned}
			&\quad \lim\limits_{\epsilon \to 0} \frac{\mathbb{P}\left(\left\|  (q, p) - (\varphi_q, \varphi_p) \right\|_{L^2_{\rho}}  \leq \epsilon\right)}{\mathbb{P}\left(\left\|  W^{\sigma} \right\|_{L^2_{\rho}}  \leq \epsilon\right)}\\
			&= {\rm exp} \left\lbrace -\frac{1}{2} \left(  \int_0^1 \left\| \sigma_q^{-1}(t) \left( \dot{\varphi}_q - \frac{\partial H}{\partial \varphi_p}(\varphi_q, \varphi_p) \right) \right\|^2_{\rho} \,{\rm d}t + \int_0^1 \left\| \sigma_p^{-1}(t) \left( \dot{\varphi}_p + \frac{\partial H}{\partial \varphi_q}(\varphi_q, \varphi_p) \right) \right\|^2_{\rho} \,{\rm d}t \right) \right\rbrace.	
		\end{aligned}\notag
	\end{displaymath}
	In general, the Onsager-Machlup functional for stochastic systems includes a correction term that accounts for path deviations induced by drift and disturbances. However, for Hamiltonian systems with conserved energy \( H \in C^3_b \), we have demonstrated that this correction term vanishes due to the symplectic structure of the system. Consequently, the Onsager-Machlup functional for Hamiltonian systems is given by:
	\begin{displaymath}
		\int_0^1 OM(\varphi, \dot{\varphi}) \,{\rm d}t = \int_0^1 \left\| \sigma_q^{-1}(t) \left( \dot{\varphi}_q - \frac{\partial H}{\partial \varphi_p}(\varphi_q, \varphi_p) \right) \right\|^2_{\rho} \,{\rm d}t + \int_0^1 \left\| \sigma_p^{-1}(t) \left( \dot{\varphi}_p + \frac{\partial H}{\partial \varphi_q}(\varphi_q, \varphi_p) \right) \right\|^2_{\rho} \,{\rm d}t.
	\end{displaymath}
\end{proof}

\section{Proof of Theorem \ref{T4.1}}

In this section, we derive the large deviation principle for stochastic Hamiltonian systems on infinite lattices by combining the Onsager-Machlup functional and Freidlin-Wentzell theory \cite{46}. The large deviation principle provides a deep insight into the probability behavior of trajectories deviating from the most probable path under stochastic perturbations, fundamentally characterizing the occurrence of rare events and their asymptotic probabilities. The rate function serves as a quantitative tool, precisely capturing the exponential decay rate of these probabilities, thus offering a detailed asymptotic representation in the limit.

\begin{proof}[Proof of Theorem \ref{T4.1}]
	Equation \(\eqref{2}\) simply specifies that the small stochastic perturbation in equation \(\eqref{1}\) is of order \( \epsilon \), where \( \epsilon \) is a small parameter. Although the forms of these two equations differ, the derivation of the Onsager-Machlup functional for equation \(\eqref{2}\) remains valid. The only distinction lies in the intensity of the stochastic noise, which is now of order \( \epsilon \) in equation \(\eqref{2}\). This scaling is reflected in the probability estimate \( \mathbb{P}\left(\left\| W^{\sigma} \right\|_{L^2_{\rho}} \leq \epsilon \right) \).
	
	Consequently, the most probable path \( \varphi(t) \) for this equation can still be derived analogously and satisfies equation \(\eqref{20}\). In this section, we focus on the large deviation behavior of the path: that is, the probability that the solution of equation \(\eqref{2}\) deviates significantly from the most probable path \( \varphi(t) \), expressed as \( \mathbb{P}\left( X(t) \in A \right) \), where \( \varphi(t) \notin A \). The set \( A \) can be decomposed into a union of tube (See Section 2.2), and therefore, we adopt definitions and techniques similar to those used in the proof of Theorem \(\ref{T3.1}\). Below, we outline the main structure of the proof.
	
	First, we define the path after applying a stochastic perturbation to \( \psi(t) \) as \( (\psi_q(t), \psi_p(t)) \), with the following expressions:
	\[
	\begin{cases}
		\Psi_q(t) = \psi_q(t) + \epsilon \int_0^t \sigma_q(s)\, {\rm d}W_q(s),\\
		\Psi_p(t) = \psi_p(t) + \epsilon \int_0^t \sigma_p(s)\, {\rm d}W_p(s).
	\end{cases}
	\]
	We then introduce a new probability measure \(\tilde{\mathbb{P}}\), under which the transformed Brownian motions are given by:
	\[
	\begin{aligned}
		\tilde{W}_q(t) &= W_q(t) - \frac{1}{\epsilon} \int_0^t \sigma_q^{-1}(s) \left( \frac{\partial H}{\partial \Psi_p}(\Psi_q, \Psi_p) - \dot{\psi}_q(s) \right)\, {\rm d}s,\\
		\tilde{W}_p(t) &= W_p(t) - \frac{1}{\epsilon} \int_0^t \sigma_p^{-1}(s) \left( -\frac{\partial H}{\partial \Psi_q}(\Psi_q, \Psi_p) - \dot{\psi}_p(s) \right)\, {\rm d}s.
	\end{aligned}
	\]
	Under the new measure \(\tilde{\mathbb{P}}\), we have
	\begin{equation*}
		\begin{cases}
			\,{\rm d} \Psi_q(t) = \frac{\partial H}{\partial \Psi_p}\left( \Psi_q, \Psi_p \right) \,{\rm d}t + \sigma_q(t) \,{\rm d} W_q(t), \\[5pt]
			\,{\rm d} \Psi_p(t) = -\frac{\partial H}{\partial \Psi_q}\left( \Psi_q, \Psi_p \right) \,{\rm d}t + \sigma_p(t) \,{\rm d} W_p(t).
		\end{cases}
	\end{equation*}
	The Radon-Nikodym derivative \(\mathcal{R} := \frac{d\tilde{\mathbb{P}}}{d\mathbb{P}}\) represents the change of measure from \(\mathbb{P}\) to \(\tilde{\mathbb{P}}\) and is given by an exponential martingale associated with the drift terms. Specifically, it takes the following form:
	\[
	\begin{aligned}
		\mathcal{R} &= \exp\left( \frac{1}{\epsilon} \left( \int_0^T \left\langle \sigma_q^{-1}(s) \left( \frac{\partial H}{\partial \Psi_p}(\Psi_q, \Psi_p) - \dot{\psi}_q(s) \right), \,{\rm d}W_q(s) \right\rangle_{\rho}\right. \right. \\
		&\qquad \qquad \left. \left. - \int_0^T \left\langle \sigma_p^{-1}(s) \left( \frac{\partial H}{\partial \Psi_q}(\Psi_q, \Psi_p) + \dot{\psi}_p(s) \right), \,{\rm d}W_p(s) \right\rangle_{\rho} \right)\right.\\ 
		&\qquad \left. - \frac{1}{2 \epsilon^2} \left( \int_0^1 \left\| \sigma_q^{-1}(s) \left( \frac{\partial H}{\partial \Psi_p}(\Psi_q, \Psi_p) - \dot{\psi}_q(s) \right) \right\|^2_{\rho}\, {\rm d}s\right. \right. \\
		&\qquad \qquad\left. \left. + \int_0^T \left\| \sigma_p^{-1}(s) \left( \frac{\partial H}{\partial \Psi_q}(\Psi_q, \Psi_p) + \dot{\psi}_p(s) \right) \right\|^2_{\rho}\, {\rm d}s \right) \right).
	\end{aligned}
	\]
	Similarly, we introduce the notation \( W^{\sigma}(t) := \left( W^{\sigma}_q(t), W^{\sigma}_p(t) \right) \), which captures the stochastic perturbation in the system:
	\[
	W^{\sigma}_q(t) := \int_0^t \sigma_q(s) \,{\rm d}W_q(s), \quad W^{\sigma}_p(t) := \int_0^t \sigma_p(s) \,{\rm d}W_p(s).
	\]
	
	Here, we calculate the probability that the trajectory of the solution \( X(t) \) of the stochastic Hamiltonian system remains sufficiently close to the reference path \( \varphi(t) \) when the noise intensity \( \epsilon \) is minimized. Specifically, we consider the probability \( \mathbb{P}(X(t) \in K(\psi, \epsilon)) \), where  
	\[  
	K(\psi, \epsilon) = \{ x - x_0 \in \mathbb{H}^1 \mid \psi - x_0 \in \mathbb{H}^1, \|x - \psi\| \leq \epsilon, \epsilon > 0 \}.  
	\]  
	By applying Girsanov's theorem, we have
	\begin{equation*}
		\begin{aligned}
			&\quad \mathbb{P}(X(t) \in K(\psi, \epsilon)) = \lim\limits_{\epsilon \to 0} \mathbb{P}\left(\left\|  (q, p) - (\psi_q, \psi_p) \right\|_{L^2_{\rho}}  \leq \epsilon \right)\\
			&= \lim\limits_{\epsilon \to 0} \tilde{\mathbb{P}}\left(\left\| (\Psi_q, \Psi_p) -(\psi_q, \psi_p) \right\|_{L^2_{\rho}}  \leq \epsilon \right)
			= \lim\limits_{\epsilon \to 0} \mathbb{E} \left( 	\mathcal{R} \mathbb{I}_{ \left\| W^{\sigma} \right\|_{L^2_{\rho}} \leq 1} \right)\\
			&= \lim\limits_{\epsilon \to 0} \mathbb{E} \left( 	\mathcal{R} \big| { \left\| W^{\sigma} \right\|_{L^2_{\rho}} \leq 1} \right) \times \mathbb{P} \left( { \left\| W^{\sigma} \right\|_{L^2_{\rho}} \leq 1 } \right) \\
			&= \lim\limits_{\epsilon \to 0} \exp\left\lbrace -\frac{1}{2\epsilon^2} \left(  \int_0^T \left\| \sigma_q^{-1}(t) \left( \dot{\psi}_q - \frac{\partial H}{\partial \psi_p}(\psi_q, \psi_p) \right) \right\|^2_{\rho}\, {\rm d}t\right. \right. \\
			&\quad \left. \left. + \int_0^T \left\| \sigma_p^{-1}(t) \left( \dot{\psi}_p + \frac{\partial H}{\partial \psi_q}(\psi_q, \psi_p) \right) \right\|^2_{\rho}\, {\rm d}t \right) \right\rbrace\\
			&\quad\times \lim\limits_{\epsilon \to 0} \mathbb{E} \left( \exp\left\lbrace \frac{1}{\epsilon^2} \sum_{i=1}^{6} B_i \right\rbrace  \big| { \left\| W^{\sigma} \right\|_{L^2_{\rho}} \leq 1 } \right)\times \mathbb{P} \left( { \left\| W^{\sigma} \right\|_{L^2_{\rho}} \leq 1 } \right),
		\end{aligned}
	\end{equation*}
	where \(B_i\) represents deviations in the path due to drift and perturbation, exhibiting stochastic characteristics. The specific form is given below:
	\begin{align*}
		B_1 &=  \epsilon \int_0^T \! \left\langle \sigma_q^{-1}(s) \frac{\partial H}{\partial \Psi_p}(\Psi_q, \Psi_p), \,{\rm d}W_q(s) \right\rangle_{\rho} \!-  \epsilon \int_0^1 \! \left\langle \sigma_p^{-1}(s)  \frac{\partial H}{\partial \Psi_q}(\Psi_q, \Psi_p), \,{\rm d}W_p(s) \right\rangle_{\rho},\\
		B_2 &= - \epsilon \int_0^T \left\langle \sigma_q^{-1}(s) \dot{\psi}_q(s), \,{\rm d}W_q(s) \right\rangle_{\rho} -  \epsilon \int_0^1 \left\langle \sigma_p^{-1}(s) \dot{\psi}_p(s), \,{\rm d}W_p(s) \right\rangle_{\rho},\\
		B_3 &= \frac{1}{2} \int_0^T \left\| \sigma_q^{-1}(s) \frac{\partial H}{\partial \psi_p}(\psi_q, \psi_p)  \right\|^2_{\rho} \,{\rm d}s - \frac{1}{2} \int_0^T \left\| \sigma_q^{-1}(s)  \frac{\partial H}{\partial \Psi_p}(\Psi_q, \Psi_p) \right\|^2_{\rho} \,{\rm d}s,\\
		B_4 &= \frac{1}{2} \int_0^T \left\| \sigma_p^{-1}(s) \frac{\partial H}{\partial \psi_q}(\psi_q, \psi_p)  \right\|^2_{\rho} \,{\rm d}s - \frac{1}{2} \int_0^T \left\| \sigma_p^{-1}(s)  \frac{\partial H}{\partial \Psi_q}(\Psi_q, \Psi_p) \right\|^2_{\rho} \,{\rm d}s,\\
		B_5 &= \int_0^T \left\langle \sigma_q^{-2}(s) \left( \frac{\partial H}{\partial \Psi_p}(\Psi_q, \Psi_p) - \frac{\partial H}{\partial \psi_p}(\psi_q, \psi_p) \right), \dot{\psi}_q(s) \right\rangle_{\rho} \,{\rm d}s,\\
		B_6 &= - \int_0^T \left\langle \sigma_p^{-2}(s) \left( \frac{\partial H}{\partial \Psi_q}(\Psi_q, \Psi_p) - \frac{\partial H}{\partial \psi_q}(\psi_q, \psi_p) \right), \dot{\psi}_p(s) \right\rangle_{\rho} \,{\rm d}s.
	\end{align*}
	We decompose the probability that the solution trajectory \(X^{\epsilon}(t)\) remains near the target path \(\psi(t)\) into three parts:
	
	1. The probability that the zero-mean Gaussian process \( W^{\sigma}(t) \) remains within a unit ball centered at the origin over the time interval \([0,T]\) is given by  
	\[
	C_{W^{\sigma}} := \mathbb{P}\left(\left\|W^{\sigma}\right\|_{L^2_{\rho}} \leq 1 \right).
	\]
	
	2. The deterministic component associated with \( \psi \) is determined by the target path itself. According to Theorem \(\ref{T3.1}\), we observe that this component corresponds to the Onsager-Machlup functional associated with Equation \(\eqref{2}\), and we denote it as
	\begin{displaymath}
		\begin{aligned}
			\exp\left\{-\frac{1}{2} \int_{0}^{T} OM(\psi, \dot{\psi}) \, \mathrm{d}s \right\}
			&:=\exp\left\{-\frac{1}{2\epsilon^2}\left(\int_{0}^{T}\left\|\sigma_q^{-1}(t)\left(\dot{\psi}_q-\frac{\partial H}{\partial \psi_p}(\psi_q,\psi_p)\right)\right\|^2_{\rho}\mathrm{d}t\right. \right. \\
			&\qquad \qquad \quad \left. \left. +\int_{0}^{T}\left\|\sigma_p^{-1}(t)\left(\dot{\psi}_p + \frac{\partial H}{\partial \psi_q}(\psi_q,\psi_p)\right)\right\|^2_{\rho}\mathrm{d}t\right)\right\}.
		\end{aligned}
	\end{displaymath}

	3. A correction part
	\[
	\mathbb{E}\left(\exp\left\{\frac{1}{\epsilon^2}\sum_{i = 1}^{6}B_i \right\} \big| {\left\|W^{\sigma}\right\|_{L^2_{\rho}} \leq 1 }\right),
	\]
	which takes into account the deviations introduced by the stochastic perturbations.

	Large deviation theory focuses on the probability of rare events occurring in a stochastic system, particularly examining how this probability decays at an exponential rate as the deviation from typical behavior increases. In this context, we are especially interested in the \(\frac{1}{\epsilon^2}\) scale. Given the boundedness of \(\sigma_{q}^{-1}(s)\) and \(\sigma_{p}^{-1}(s)\), along with the fact that \(H \in C^3_b\big(\ell^2_{\rho}(\mathbb{Z}^m; M); \mathbb{R}\big)\), we can infer that
	\begin{equation}\label{24}
		\limsup\limits_{\substack{\epsilon \to 0 \\ \delta \to 0}} \mathbb{E}\left({\rm exp}\left\{ cB_{1} \right\} \big| { \left\| W^{\sigma} \right\|_{L^2_{\rho}} \leq 1} \right) = 1,
	\end{equation}
	and
	\begin{equation}\label{25}
		\limsup\limits_{\substack{\epsilon \to 0 \\ \delta \to 0}} \mathbb{E}\left({\rm exp}\left\{ cB_{2} \right\} \big| { \left\| W^{\sigma} \right\|_{L^2_{\rho}} \leq 1 } \right) = 1
	\end{equation}
	for all $c \in \mathbb{R}$.
	
	For the third term $B_3$, 
	\begin{displaymath}
		\begin{aligned}
			B_3 &= \frac{1}{2} \int_0^T \left\| \sigma_q^{-1}(s) \frac{\partial H}{\partial \psi_p}(\psi_q, \psi_p)  \right\|^2_{\rho} \,{\rm d}s - \frac{1}{2} \int_0^1 \left\| \sigma_q^{-1}(s)  \frac{\partial H}{\partial \Psi_p}(\Psi_q, \Psi_p) \right\|^2_{\rho} \,{\rm d}s
			\\ &\leq \frac{1}{2} \int_{0}^{T} {\sigma_q^{-2}(s) \left\|  \frac{\partial H}{\partial \psi_p}(\psi_q, \psi_p) - \frac{\partial H}{\partial \Psi_p}(\Psi_q, \Psi_p) \right\|^2_{\rho} }\,{\rm d}s\\
			&\quad + \int_{0}^{T} {\sigma_q^{-2}(s) \left\|\frac{\partial H}{\partial \psi_p}(\psi_q, \psi_p) - \frac{\partial H}{\partial \Psi_p}(\Psi_q, \Psi_p) \right\|_{\rho} \left\| \frac{\partial H}{\partial \Psi_p}(\Psi_q, \Psi_p) \right\|_{\rho}} \,{\rm d}s.
		\end{aligned}
	\end{displaymath}
	Using that $\frac{\partial H}{\partial p}$ is Lipschitz continuous, we have
	\begin{equation}\label{26}
		\begin{aligned}
			&\quad \left\| \frac{\partial H}{\partial \Psi_p}(\Psi_q, \Psi_p)  -  \frac{\partial H}{\partial \psi_p}(\psi_q, \psi_p)\right\|_{\rho}\\
			&= \left\| \frac{\partial H}{\partial \left( \psi_p + \epsilon W^{\sigma}_p \right)}((\psi_q + \epsilon W^{\sigma}_q), (\psi_p + \epsilon W^{\sigma}_p))-  \frac{\partial H}{\partial \psi_p}(\psi_q, \psi_p)\right\|_{\rho}\\
			&\leq L \epsilon \left\| W^{\sigma} \right\|_{\rho}.
		\end{aligned}
	\end{equation}
	Inequality $\eqref{26}$ and the boundedness of $\frac{\partial H}{\partial y_p}(y_q, y_p)$ and $\sigma_q^{-1}(t)$ imply that
	\begin{equation}\label{27}
		\limsup\limits_{\substack{\epsilon \to 0 \\ \delta \to 0}} \mathbb{E}\left({\rm exp}\left\{ cB_3 \right\} \big| { \left\| W^{\sigma} \right\|_{L^2_{\rho}} \leq 1} \right) = 1
	\end{equation}
	for all $c \in \mathbb{R}$.
	
	For the fourth term $B_4$, employing the same proof technique as for the third term $B_3$, we have
	\begin{equation}\label{28}
		\limsup\limits_{\substack{\epsilon \to 0 \\ \delta \to 0}} \mathbb{E}\left({\rm exp}\left\{ cB_4 \right\} \big| { \left\| W^{\sigma} \right\|_{L^2_{\rho}} \leq 1} \right) = 1
	\end{equation}
	for all $c \in \mathbb{R}$.
	
	For the fifth term $B_5$, applying inequality $\eqref{26}$ and the boundedness of $\dot{\varphi}_q(t)$ and $\sigma_q^{-1}(t)$, we have
	\begin{displaymath}
		\begin{aligned}
			B_5 &= \int_0^T \left\langle \sigma_q^{-2}(s) \left( \frac{\partial H}{\partial \Psi_p}(\Psi_q, \Psi_p) - \frac{\partial H}{\partial \psi_p}(\psi_q, \psi_p) \right), \dot{\psi}_q(s) \right\rangle_{\rho} \,{\rm d}s\\
			& \leq C \int_0^T \left\|  \frac{\partial H}{\partial \Psi_p}(\Psi_q, \Psi_p) - \frac{\partial H}{\partial \psi_p}(\psi_q, \psi_p) \right\|_{\rho} \,{\rm d}s\\
			& \leq C\int_0^T L \left\| W^{\sigma} \right\|_{\rho} \,{\rm d}s\\
			& \leq CLT^{\frac{1}{2}}\left\| W^{\sigma} \right\|_{L^2_{\rho}}.
		\end{aligned}
	\end{displaymath}
	Thus,
	\begin{equation}\label{29}
		\limsup\limits_{\substack{\epsilon \to 0 \\ \delta \to 0}} \mathbb{E}\left({\rm exp}\left\{ cB_5 \right\} \big| { \left\| W^{\sigma} \right\|_{L^2_{\rho}} \leq 1} \right) = 1
	\end{equation}
	for all $c \in \mathbb{R}$.
	
	For the sixth term $B_6$, employing the same proof technique as for the fifth term $B_5$, we have
	\begin{equation}\label{30}
		\limsup\limits_{\substack{\epsilon \to 0 \\ \delta \to 0}} \mathbb{E}\left({\rm exp}\left\{ cB_6 \right\} \big| { \left\| W^{\sigma} \right\|_{L^2_{\rho}} \leq 1} \right) = 1
	\end{equation}
	for all $c \in \mathbb{R}$.
	
	Building on the results of Equations \eqref{24}, \eqref{25}, and \eqref{27}–\eqref{30}, we derive the following asymptotic approximation:  
	\begin{equation}
		\begin{aligned}
			\epsilon^2 \ln \mathbb{P}(X(t) \in K(\psi, \epsilon)) &= \epsilon^2 \ln \left( \exp\left\{-\frac{1}{2} \int_{0}^{T} OM(\psi, \dot{\psi}) \, \mathrm{d}s \right\}\right. \\
			&\qquad \left. \times \mathbb{E}\left(\exp\left\{\frac{1}{\epsilon^2}\sum_{i = 1}^{6}B_i \right\} \big| {\left\|W^{\sigma}\right\|_{L^2_{\rho}} \leq 1 }\right) \times C_{W^{\sigma}} \right),\\
			&= \epsilon^2 \ln \left( \exp\left\{-\frac{1}{2} \int_{0}^{T} OM(\psi, \dot{\psi}) \, \mathrm{d}s \right\}\right) \\
			&\quad + \epsilon^2 \ln \left( \mathbb{E}\left(\exp\left\{\frac{1}{\epsilon^2}\sum_{i = 1}^{6}B_i \right\} \big| {\left\|W^{\sigma}\right\|_{L^2_{\rho}} \leq 1 }\right)\right)\\
			&\quad + \epsilon^2 \ln \left( C_{W^{\sigma}} \right),\\
			&\approx \epsilon^2 \ln \left(  \exp\left\{-\frac{1}{2} \int_{0}^{T} OM(\psi, \dot{\psi}) \, \mathrm{d}s \right\} \right).
		\end{aligned} 
	\end{equation}
	This is because the latter two terms vanish as $\epsilon$ approaches zero. Proceeding from this foundation, for an arbitrary measurable set \( A \), we compute \( \epsilon^2 \ln \mathbb{P}(X(t) \in A) \). Specifically, for any continuous function \( \varphi \in A \), we focus exclusively on the case where \( \varphi - x_0 \in \mathbb{H}^{1} \); for all other cases, we set \( J = \infty \) by definition.
	
	Upper Bound Estimation. Let \( F \subset \mathbb{H}^1 \) be a closed set. By the continuity of \( U(\psi) \), there exists a sequence of compact sets \( \{K_n\} \) satisfying \( K_n \subset F \) with:
	\[
	\inf_{\psi \in K_n} U(\psi) \to \inf_{\psi \in F} U(\psi) \quad \text{as } n \to \infty.
	\]
	For each \( K_n \) and \( \epsilon > 0 \), we construct a finite covering:
	\[
	K_n \subset \bigcup_{j=1}^N K(\psi_j, \epsilon), \quad \{\psi_j\}_{j=1}^N \subset K_n.
	\]
	Applying the union bound:
	\[
	\mathbb{P}(X(t) \in K_n) \leq \sum_{j=1}^N \mathbb{P}(X(t) \in K(\psi_j, \epsilon)).
	\]
	Taking logarithmic scaling:
	\begin{align*}
		\epsilon^2 \ln \mathbb{P}(X(t) \in K_n) &\leq \epsilon^2 \ln \left( N \cdot \max_{1 \leq j \leq N} \mathbb{P}(X(t) \in K(\psi_j, \epsilon)) \right) \nonumber \\
		&\leq \epsilon^2 \ln N + \max_{1 \leq j \leq N} \epsilon^2 \ln \mathbb{P}(X(t) \in K(\psi_j, \epsilon)).
	\end{align*}
	Taking the limsup as \( \epsilon \to 0 \):
	\begin{equation*}
		\limsup_{\epsilon \to 0} \epsilon^2 \ln \mathbb{P}(X(t) \in K_n) \leq -\min_{1 \leq j \leq N} \frac{1}{2} \int_0^T U(\psi_j) \,{\rm d}s.
	\end{equation*}
	For any \( \delta > 0 \), select \( K_n \subset F \) such that:
	\[
	\inf_{\psi \in F} \int_0^T U(\psi) \,{\rm d}s + \delta \geq \inf_{\psi \in K_n} \int_0^T U(\psi) \,{\rm d}s \geq \inf_{\psi \in F} \int_0^T U(\psi) \,{\rm d}s - \delta.
	\]
	By combining this with the exponential decay of probabilities outside the set \( K_n \), it follows that:
	\begin{equation}\label{eq:upper_bound}
		\limsup_{\epsilon \to 0} \epsilon^2 \ln \mathbb{P}(X(t) \in F) \leq -\inf_{\psi \in F} \frac{1}{2} \int_0^T U(\psi) \,{\rm d}s.
	\end{equation}
	
	Lower Bound Estimation. Let \( G \subset \mathbb{R}^n \) be open. For any \( \psi \in G \):
	\[
	\liminf_{\epsilon \to 0} \epsilon^2 \ln \mathbb{P}(X(t) \in G) \geq \liminf_{\epsilon \to 0} \epsilon^2 \ln \mathbb{P}(X^\epsilon \in K(\psi, \epsilon)).
	\]
	Taking supremum over \( \psi \in G \):
	\begin{equation}\label{eq:lower_bound}
		\liminf_{\epsilon \to 0} \epsilon^2 \ln \mathbb{P}(X(t) \in G) \geq -\inf_{\psi \in G} \frac{1}{2} \int_0^T U(\psi) \,{\rm d}s.
	\end{equation}
	
	Combining \eqref{eq:upper_bound} and \eqref{eq:lower_bound}, the rate function is:
	\begin{align}
		J(\psi) &= \frac{1}{2} \bigg( \int_0^T \left\| \sigma_q^{-1}(t) \left( \dot{\psi}_q - \frac{\partial H}{\partial \psi_p}(\psi_q, \psi_p) \right) \right\|^2 \,{\rm d}t \nonumber \\
		&\quad + \int_0^T \left\| \sigma_p^{-1}(t) \left( \dot{\psi}_p + \frac{\partial H}{\partial \psi_q}(\psi_q, \psi_p) \right) \right\|^2 \,{\rm d}t \bigg).
	\end{align}
\end{proof}

\section{ Proof of Theorem \ref{T5.1}}

In this section, we integrate the Onsager-Machlup functional, the most probable path, and the large deviation principle with the KAM theory proposed in \cite{25} to investigate the persistence of invariant tori in the stochastic nonlinear Schr{\"{o}}dinger equation on infinite lattices. Within the framework of the most probable path, we analyze the impact of stochastic perturbations on the system's stability and trajectory evolution, and quantify the probability of invariant tori persistence using the large deviation principle. This approach provides novel insights and tools for understanding invariant structures in physical systems.

\begin{proof}[Proof of Theorem \ref{T5.1}]
	We consider the nonlinear Schrödinger equation \eqref{50} as a Hamiltonian system formulated on a suitable phase space $\mathcal{P}$. For instance, we may take \( \mathcal{P} = W_{0}^{1}([0, \pi]) \), the Sobolev space of all complex-valued \( L^{2} \)-functions on \([0, \pi]\) with an \( L^{2} \)-derivative and vanishing boundary values. Equipped with the inner product
	\[
	\langle u, v \rangle = \mathrm{Re} \int_{0}^{\pi} u \overline{v} \, \mathrm{d}x.
	\]
	The Hamiltonian of the nonlinear Schr{\"{o}}dinger equation \eqref{50} is given by
	\[
	H = \frac{1}{2} \langle A u, u \rangle + \frac{1}{2} \int_{0}^{\pi} g(|u|^2) \, \mathrm{d}x,
	\]
	where \( A = -\frac{d^2}{dx^2} + m \) is the operator acting on the wave function, and \( g(s) = \int_{0}^{s} f(z) \, \mathrm{d}z \) represents the nonlinear interaction. To express this Hamiltonian in terms of infinitely many coordinates, we introduce the representation:
	\[
	u = \mathcal{S} q = \sum_{j=1}^{\infty} q_j \phi_j, \quad \phi_j = \sqrt{\frac{2}{\pi}} \sin(jx), \quad j \geq 1.
	\]
	Here, the coordinates \( q = (q_1, q_2, \ldots) \) belong to the Hilbert space \( \ell^{a,p} \), consisting of all complex-valued sequences satisfying
	\[
	\|q\|_{a,p}^2 = \sum_{j=1}^{\infty} |q_j|^2 j^{2p} e^{2j a} < \infty,
	\]
	for \( a > 0 \) and \( p \geq \frac{1}{2} \) fixed later. In these coordinates, the Hamiltonian becomes
	\[
	H = \Lambda + G = \frac{1}{2} \sum_{j=1}^{\infty} \lambda_j |q_j|^2 + \frac{1}{2} \int_0^\pi g(|\mathcal{S} q|^2) \, \mathrm{d}x.
	\]
	The phase space \( \ell^{a,p} \) is equipped with the symplectic structure
	\[
	\frac{i}{2} \sum_j dq_j \wedge d\overline{q_j}.
	\]
	Hamilton's equations of motion are then given by
	\[
	\dot{q}_j = 2i \frac{\partial H}{\partial \overline{q_j}}, \quad j \geq 1,
	\]
	which represent the classical Hamiltonian dynamics. Expressing \( q_j \) and $\overline{q}_j$ in terms of their real and imaginary parts,
	\[
	q_j = x_j + i y_j, \quad \overline{q}_j = x_j - i y_j,
	\]
	where \( x_j, y_j \in \mathbb{R} \), the Hamiltonian \( H \) can be written as \( H = H(x_1, x_2, \ldots, y_1, y_2, \ldots) \).
	
	Computing the partial derivatives using the chain rule, we obtain:
	\[
	\frac{\partial H}{\partial \overline{q}_j} = \frac{1}{2} \frac{\partial H}{\partial x_j} - \frac{1}{2i} \frac{\partial H}{\partial y_j},
	\]
	which leads to the equations of motion:
	\[
	\dot{q}_j = i \frac{\partial H}{\partial x_j} + \frac{\partial H}{\partial y_j}.
	\]
	Separating real and imaginary parts, we obtain the system:
	\[
	\begin{cases}
		\dot{x}_j = \frac{\partial H}{\partial y_j}, \\
		\dot{y}_j = -\frac{\partial H}{\partial x_j},
	\end{cases} \quad j \geq 1.
	\]
	In vector form, for the infinite-dimensional vectors \( X = (x_1, x_2, \ldots) \) and \( Y = (y_1, y_2, \ldots) \), the equations of motion are:
	\begin{equation}\label{45}
		\begin{cases}
			\dot{X} = \frac{\partial H}{\partial Y},\\
			\dot{Y} = -\frac{\partial H}{\partial X}.
		\end{cases}
	\end{equation}
	
	Introducing stochastic perturbations into equation \eqref{50}, we obtain the stochastic nonlinear Schr{\"{o}}dinger equation:
	\begin{equation}\label{46}
		i \frac{\partial u(x,t)}{\partial t} = \frac{\partial^2 u(x,t)}{\partial x^2} - m u(x,t) - f(|u(x,t)|^2) u(x,t) + \epsilon \dot{\eta}(t,x), \quad t \in [0,T],~ x \in [0,\pi],
	\end{equation}
	where \( \epsilon \) denotes the noise intensity, and \( \eta(t, x) \) is a stochastic process defined on a probability space \((\Omega, \mathcal{F}, P)\). We assume that \( \eta(t, x) \) has the expansion
	\[
	\eta(t, x) = \sum_{j=1}^{\infty} \sigma_j(t) W_j(t) \phi_j(x),
	\]
	where \( \sigma_j(t) = \sigma_j^R(t) + i \sigma_j^I(t) \) are complex-valued functions, and \( W_j(t) \) are independent standard Brownian motions for each \( j \geq 1 \). This allows us to reformulate the stochastic system as an infinite-dimensional vector system:
	\begin{equation}\label{48}
		\begin{cases}
			\dot{X} = \frac{\partial H}{\partial Y} + \epsilon \sigma_R(t) \dot{W}(t),\\
			\dot{Y} = -\frac{\partial H}{\partial X} + \epsilon \sigma_I(t) \dot{W}(t),
		\end{cases}
	\end{equation}
	where \( X = (x_1, x_2, \ldots) \) and \( Y = (y_1, y_2, \ldots) \) are infinite-dimensional vectors, \( \sigma_R(t) \) and \( \sigma_I(t) \) are diagonal matrices of noise coefficients, and \( W(t) = (W_j(t))_{j \geq 1} \) is an infinite-dimensional vector of independent standard Brownian motions:
	\[
	\sigma_R(t) = \operatorname{diag}(\sigma^R_j(t))_{j \geq 1}, \quad 
	\sigma_I(t) = \operatorname{diag}(\sigma^I_j(t))_{j \geq 1}, \quad
	W(t) = (W_j(t))_{j \geq 1}.
	\]
	
	According to Theorem \ref{T3.1}, when equation \eqref{48} satisfies conditions \( (C1) \) and \( (C2) \), its Onsager-Machlup functional is given by:
	\[
	\begin{aligned}
		\int_0^{T} OM(\varphi, \dot{\varphi}) \, \mathrm{d}t &= \frac{1}{\epsilon^2} \left( \int_0^{T} \left\| \sigma_R^{-1}(s) \left( \dot{\varphi}_X(t) - \frac{\partial H}{\partial \varphi_Y}(\varphi_X, \varphi_Y) \right) \right\|^2_{a,p} \, \mathrm{d}t \right. \\
		&\quad + \left. \int_0^{T} \left\| \sigma_I^{-1}(s) \left( \dot{\varphi}_Y(t) + \frac{\partial H}{\partial \varphi_X}(\varphi_X, \varphi_Y) \right) \right\|^2_{a,p} \, \mathrm{d}t \right),
	\end{aligned}
	\]
	where \( \varphi = (\varphi_X, \varphi_Y) \) represents the state variables, \( H(\varphi_X, \varphi_Y) \) is the Hamiltonian, and \( \sigma_R^{-1}(s) \), \( \sigma_I^{-1}(s) \) are the inverses of the noise intensity matrices.
	
	By minimizing the Onsager-Machlup functional, the most probable transition path \( \varphi = (\varphi_X, \varphi_Y) \) for the stochastic system \eqref{48} satisfies the corresponding deterministic system \eqref{45} without stochastic perturbations. This implies that the real and imaginary parts of the solution to the stochastic nonlinear Schr{\"{o}}dinger equation \eqref{46} are determined by the deterministic equation \eqref{50}.
	
	Furthermore, Theorem \ref{T4.1} establishes the large deviation principle in the following form: for any measurable rare event set \( A \), the probability that the solution \( q(t) = (q_n(t))_{n \in \mathbb{Z}} \) of system \eqref{48} lies in \( A \) satisfies
	\[
	\epsilon^2 \ln \mathbb{P}(q(t) \in A ) \approx - \inf_{\psi \in A}  J(\psi),
	\]
	where the \textit{rate function} \( J(\psi) \) is given by
	\begin{align*}
		J(\psi) &= \frac{1}{2} \int_0^{T} \left\| \sigma_R^{-1}(s) \left( \dot{\varphi}_X(t) - \frac{\partial H}{\partial \varphi_Y}(\varphi_X, \varphi_Y) \right) \right\|^2_{a,p} \, \mathrm{d}t
		+ \frac{1}{2} \int_0^{T} \left\| \sigma_I^{-1}(s) \left( \dot{\varphi}_Y(t) + \frac{\partial H}{\partial \varphi_X}(\varphi_X, \varphi_Y) \right) \right\|^2_{a,p} \, \mathrm{d}t.
	\end{align*}
	Alternatively, we may reformulate it as a stochastic nonlinear Schrödinger equation. In this setting, the probability that the stochastic trajectory $u^{*}(t,x)$ of equation \eqref{50} deviates from the invariant torus is governed by the following large deviation principle:
	\begin{equation*}
		\epsilon^2 \ln \mathbb{P}(u^{*}(t,x) \in D ) \approx - \inf_{\psi \in D}  J(\psi),
	\end{equation*}
	where \( \psi \in D\) denotes an arbitrarily continuous function, and \( D \subset M \) denote an arbitrary measurable set. Furthermore, the rate function \( J(\psi) \) is given by:
	\begin{equation*}
		\begin{aligned}
			J(\psi)=
			\begin{cases}
				\frac{1}{2} \int_0^T \int_0^{\pi} \left\| \Sigma^{-1}(t) \left(  \frac{\partial \psi(x,t)}{\partial t} + i \left( \frac{\partial^2 \psi(x,t)}{\partial x^2} - m \psi(x,t) - f(|\psi(x,t)|^2) \psi(x,t)\right)  \right) \right\|^2_{a,p}\,{\rm d}x \,{\rm d}t , & \text{if } \psi - x_0 \in \mathcal{H}^1;\\
				+\infty, & \text{otherwise}.
			\end{cases}
		\end{aligned}
	\end{equation*}
	where
	$
	\Sigma^{-1}(t) v := \sum_{j=1}^{\infty} \left( \frac{\langle \operatorname{Re} v, \phi_j \rangle}{\sigma_j^R(t)} \;+\; i \,\frac{\langle \operatorname{Im} v, \phi_j \rangle}{\sigma_j^I(t)} \right) \phi_j(x)
	$
	is defined as the inverse of the noise operator, with $v(x)$ being a complex-valued function. And $\mathcal{H}^1 = H^1(0,T;\ell^{a,p}(0,\pi)) $ is a Bochner–Sobolev space, meaning that it is a Sobolev space in the time variable whose values lie in the spatial space $\ell^{a,p}(0,\pi)$.
	
	In summary, we have derived the Onsager-Machlup functional, the most probable transition path, and the large deviation principle for the stochastic nonlinear Schr{\"{o}}dinger equation on infinite lattices. This framework allows us to analyze the properties of solution trajectories within a probabilistic context, particularly focusing on the preservation of low-dimensional invariant tori. To achieve this, we leverage the KAM theory for the deterministic nonlinear Schr{\"{o}}dinger equation, a well-studied problem with significant results established by various researchers. A similar KAM approach is adopted in our analysis (see, for example, \cite{74}).
	
	Firstly, we examine the nonlinear term \( |u|^2 u \) in the stochastic nonlinear Schr{\"{o}}dinger equation, which determines the coefficients in the Birkhoff normal form. We find that
	\[
	G = \frac{1}{4} \int_{0}^{\pi} |u(x)|^4 \, \mathrm{d}x = \frac{1}{4} \sum_{i,j,k,l} G_{ijkl} q_i q_j \overline{q_k} \overline{q_l},
	\]
	where
	\[
	G_{ijkl} = \int_{0}^{\pi} \phi_i \phi_j \phi_k \phi_l \, \mathrm{d}x.
	\]
	It is straightforward to verify that \( G_{ijkl} = 0 \) unless there exists a combination of signs such that \( i \pm j \pm k \pm l = 0 \).
	
	Therefore, for the Hamiltonian \( H = \Lambda + G \), there exists a real analytic, symplectic change of coordinates \( \Gamma \) in a neighborhood of the origin in \( \ell^{a,p} \), which transforms \( H \) into its Birkhoff normal form up to fourth order for all real values of \( m \). Specifically,
	\[
	H \circ \Gamma = \Lambda + \overline{G} + K,
	\]
	where \( X_{\overline{G}} \) and \( X_{K} \) are real analytic vector fields in a neighborhood of the origin in \( \ell^{a,p} \), 
	\[
	\overline{G} = \frac{1}{2} \sum_{i,j \geq 1} \overline{G}_{ij} |q_i|^{2} |q_j|^{2}, \quad |K| = O\left( \|q\|^{6}_{\ell^{a,p}} \right),
	\]
	and the coefficients \( \overline{G}_{ij} \) are uniquely determined as \( \overline{G}_{ij} = \frac{4 - \delta_{ij}}{4\pi} \).
	
	In complex coordinates \( q = (\hat{q}, \check{q}) \) on \( \ell^{a,p} \), where \( \hat{q} = (q_1, \ldots, q_n) \) and \( \check{q} = (q_{n+1}, q_{n+2}, \ldots) \), and with
	\[
	I = \frac{1}{2} (|q_1|^2, \ldots, |q_n|^2), \quad Z = \frac{1}{2} (|q_{n+1}|^2, |q_{n+2}|^2, \ldots),
	\]
	the normal form consists of the terms
	\[
	\Lambda = \langle \alpha, I \rangle + \langle \beta, Z \rangle, \quad Q = \frac{1}{2} \langle A I, I \rangle + \langle B I, Z \rangle,
	\]
	where \( \alpha \) and \( \beta \) are vectors with constant coefficients, and \( A \) and \( B \) are matrices with constant coefficients. The equations of motion are:
	\[
	\dot{\hat{q}} = i \, \text{diag}(\alpha + A I + B^T Z)\hat{q}, \quad \dot{\check{q}} = i \, \text{diag}(\beta + B I)\check{q}.
	\]
	These equations possess a complex \( n \)-dimensional invariant manifold \( E = \{\check{q} = 0\} \), which is completely filled, up to the origin, by the invariant tori
	\[
	T(I) = \left\{ \hat{q} : |\hat{q}_j|^2 = 2 I_j \quad \text{for} \quad 1 \leq j \leq n \right\}, \quad I \in \mathbb{R}^n.
	\]
	On \( T(I) \) and in its normal space, the flows are given by
	\begin{equation*}
		\begin{cases}
			\dot{\hat{q}} = i \, \text{diag}(\omega(I)) \hat{q}, \quad \omega(I) = \alpha + A I,\\
			\dot{\check{q}} = i \, \text{diag}(\Omega(I)) \check{q}, \quad \Omega(I) = \beta + B I.
		\end{cases}
	\end{equation*}
	These equations are linear and diagonal. In particular, since \( \Omega(I) \) is real, \( \check{q} = 0 \) is an elliptic fixed point, all tori are linearly stable, and all their orbits have zero Lyapunov exponents. We refer to \( T(I) \) as an elliptic rotational torus with frequencies \( \omega(I) \).
	
	Due to resonance, the manifold \( E \) containing higher-order terms generally does not exist. Instead, we aim to show that a significant portion of \( E \) can form an invariant Cantor manifold \( \mathcal{E} \). Specifically, there exists a family of \( n \)-dimensional tori
	\[
	T[C] = \bigcup_{I \in C} T(I) \subset E,
	\]
	where \( C \) is a Cantor set in \( \mathbb{R}^n \), and there exists a Lipschitz continuous embedding
	\[
	\varPsi : T[C] \to \ell^{a,p},
	\]
	such that the restriction of \( \varPsi \) to each torus \( T(I) \) is an embedding of an elliptic rotational \( n \)-torus associated with the Hamiltonian \( H \). The image \( \mathcal{E} \) of \( T[C] \) under \( \varPsi \) is called the Cantor manifold of elliptic rotational \( n \)-tori provided by the embedding \( \varPsi : T[C] \to \mathcal{E} \). Moreover, the Cantor set \( C \) has full density at the origin, the embedding \( \varPsi \) is close to the inclusion map \( \varPsi_0 : E \to \ell^{a,p} \), and the Cantor manifold \( \mathcal{E} \) is tangent to \( E \) at the origin. The existence of \( \mathcal{E} \) is established under the following assumptions:
	\begin{itemize}
		\item[(C3)] Nondegeneracy. The normal form \( \Lambda + Q \) is nondegenerate in the sense that
		\begin{equation*}
			\det A \neq 0, \qquad \langle l, \beta \rangle \neq 0, \qquad
			\langle k, \omega(I) \rangle + \langle l, \Omega(I) \rangle \neq 0,
		\end{equation*}
		for all \( (k, l) \in \mathbb{Z}^n \times \mathbb{Z}^{\infty} \) with \( 1 \leq \left| l\right| \leq 2 \), where \( \omega = \alpha + AI \) and \( \Omega = \beta + BI \).
		
		\item[(C4)] Spectral Asymptotics. There exist \( d \geq 1 \) and \( \delta < d - 1 \) such that
		\[
		\beta_j = j^d + \cdots + O(j^{\delta}),
		\]
		where the dots denotes terms of order less than \( d \) in \( j \).
		
		\item[(C5)] Regularity.
		\begin{align*}
			X_Q, X_R \in A(\ell^{a,p}, \ell^{a,\bar{p}}), \qquad
			\begin{cases}
				\bar{p} \geq p & \text{for } d > 1, \\
				\bar{p} > p & \text{for } d = 1,
			\end{cases}
		\end{align*}
		where $A(\ell^{a,p}, \ell^{a,\bar{p}})$ denotes the class of all maps from some neighbourhood of the origin in $\ell^{a,p}$ into $\ell^{a,\bar{p}}$, which are real analytic in the real and imaginary parts of the complex coordinate $q$.
	\end{itemize}
	
	Building on Assumptions (C3)-(C5), we now describe the structure of the solutions. Consider the Hamiltonian $H = \Lambda + Q + R$. If the remainder term satisfies
	$$
	|R| = O\!\left(\|q\|^{q}_{\ell^{a,p}}\right) + O\!\left(\|\hat{q}\|^{4}_{\ell^{a,\hat{p}}}\right),
	$$
	together with the condition
	$$
	q > 4 + \frac{4 - \Delta}{\kappa}, \qquad \Delta = \min(\bar{p} - p, 1),
	$$
	then one can establish the existence of a Cantor-type manifold $\mathcal{E}$ consisting of real-analytic, elliptic Diophantine $n$-tori. More precisely, $\mathcal{E}$ can be represented by a Lipschitz continuous embedding
	$$
	\varPsi : T[C] \to \mathcal{E},
	$$
	where the set $C$ has full density at the origin. The embedding $\varPsi$ is close to the natural inclusion map $\varPsi_{0}$, in the sense that
	$$
	\|\varPsi - \varPsi_{0}\|_{\ell^{a,p}, B^{-1}\Gamma[C]} = O\!\left(r^{\sigma}\right), 
	\qquad 
	\sigma = \frac{q}{2} - \frac{\kappa + 1 - \Delta/4}{\kappa} > 1.
	$$
	As a consequence, the manifold $\mathcal{E}$ is tangent to $E$ at the origin. For a more detailed proof, we refer the reader to the Cantor Manifold Theorem in \cite{74}.

	In summary, we have established the KAM theorem stated in Theorem \ref{T2.6}. On the basis of this result, we have completed the proof of the stochastic version of the KAM theorem, thereby confirming the persistence of invariant tori for the one-dimensional stochastic nonlinear Schrödinger equation. This concludes the proof of our main result.
	
\end{proof}

﻿

\bibliographystyle{plain}
\bibliography{Ref}

\end{document}